%% file: arxiv.tex
\documentclass[conference]{IEEEtran}
\ifCLASSINFOpdf

\else

\fi

\usepackage{epsfig} 
\usepackage{amsmath,amssymb, amsxtra, mathrsfs}
\usepackage{xcolor}
\usepackage{multirow}
\usepackage{caption}
\usepackage{subcaption}
\usepackage{url}
\usepackage{booktabs,tabularray,tabu}
\usepackage{breqn}

\usepackage{algorithm, algpseudocode}

\usepackage{pgfplots,tikz}
\pgfplotsset{compat=1.15}
\usepackage{tikz}
\usetikzlibrary{shapes,arrows,spy,positioning,patterns,shapes.geometric}
\usetikzlibrary{external}
\newtheorem{definition}{Definition}
\newtheorem{definition*}{Definition}

\begin{document}
\pagenumbering{Arabic}
\title{Water Demand Maximization: Quick Recovery of Nonlinear Physics Solutions}


\author{\IEEEauthorblockN{Sai Krishna K. Hari}
\IEEEauthorblockA{Applied Mathematics \& Plasma Physics (T-5)\\
Los Alamos National Laboratory\\
Los Alamos, NM - 87544\\
Email: hskkanth@lanl.gov}
\and
\IEEEauthorblockN{Russell Bent}
\IEEEauthorblockA{Applied Mathematics \& Plasma Physics (T-5)\\
Los Alamos National Laboratory\\
Los Alamos, NM - 87544\\
Email: rbent@lanl.gov}
}


%


\maketitle

\begin{abstract}
Determining the maximum demand a water distribution network can satisfy is crucial for ensuring reliable supply and planning network expansion.
This problem, typically formulated as a mixed-integer nonlinear program (MINLP), is computationally challenging. A common strategy to address this challenge is to solve mixed-integer linear program (MILP) relaxations derived by partitioning variable domains and constructing linear over- and under-estimators to nonlinear constraints over each partition.
While MILP relaxations are easier to solve up to a modest level of partitioning, their solutions often violate nonlinear water flow physics.
Thus, recovering feasible MINLP solutions from the MILP relaxations is crucial for enhancing MILP-based approaches.
In this paper, we propose a robust solution recovery method that efficiently computes feasible MINLP solutions from MILP relaxations, regardless of partition granularity. Combined with iterative partition refinement, our method generates a sequence of feasible solutions that progressively approach the optimum.
Through extensive numerical experiments, we demonstrate that our method outperforms baseline methods and direct MINLP solves by consistently recovering high-quality feasible solutions with significantly reduced computation times.

\end{abstract}

\IEEEpeerreviewmaketitle

\section{Introduction}
Reliable and efficient operation of water distribution networks (WDNs) is vital for ensuring adequate supply, identifying system vulnerabilities, and planning expansions that enhance system resilience \cite{yazdani2011resilience}. Within this context, a key planning problem is determining the maximum demand that a WDN can satisfy under hydraulic and operational limits, as its solution enables operators to evaluate system capacity, 
anticipate future constraints, and coordinate with interdependent infrastructures such as power and transportation networks \cite{hari2024joint}.
The demand maximization problem is often formulated as a mixed-integer nonlinear program (MINLP), which captures the nonlinear, nonconvex water flow physics and the discrete on/off operation of pumps and valves.
These MINLPs are computationally challenging to solve, especially for large-scale networks.

A common strategy adopted in the literature to improve tractability is to construct linear or piecewise-linear relaxations of the nonlinear models \cite{tasseff2024polyhedral, hari2023relaxations, bonvin2021pump}. These relaxations typically partition variable domains by introducing disjunctive variables and constructing linear over- and under-estimators for the nonlinear constraints within the partitioned space. This approach yields mixed-integer linear programming (MILP) relaxations of the original MINLPs.

These MILP relaxations are relatively easier to solve up to a modest level of partitioning, but tractability diminishes rapidly as the partitions are refined. However, their solutions often violate nonlinear hydraulic constraints, rendering them infeasible for the original MINLP. 
Although the relaxations provide valuable bounds that can be embedded in global optimization algorithms (such as spatial branch-and-bound) to enhance efficiency, in many applications, including water network planning, a quick feasible solution is often more valuable than a computationally intensive global optimum \cite{menke2016exploring}.
Consequently, efficiently recovering high-quality feasible MINLP solutions from MILP relaxations has become an important step toward enhancing the appeal of MILP-based methods but remains an underexplored challenge.
Baseline recovery methods \cite{borraz2016convex} typically fix the integer variables of the MINLP to the corresponding solution values of their MILP relaxations and then solve the resulting nonlinear program (NLP). These methods often fail when the NLP is infeasible, particularly under coarse relaxations.

To overcome the challenges of infeasibility and intractability, this paper develops a robust recovery algorithm that yields feasible MINLP solutions from a MILP with an arbitrary (even coarse) level of partitioning.
We base the recovery step on the observation that while relaxations may produce \text{infeasible} integer solutions, they still contain valuable information. Often, a feasible integer solution is not far away, and re-solving the MILP with a granular partitioning is not necessary to find such a feasible solution.
We combine this observation with two parallel steps: neighborhood search and iterative partition refinement, to develop a \textit{robust recovery algorithm} that quickly computes high-quality feasible solutions. The algorithm exploits the tradeoffs between relaxation strength and computational tractability, and
when combined with a sequence of converging piecewise relaxations \cite{sundar2022sequence}, produces a sequence of feasible solutions that progressively approach the optimal solution.
The effectiveness of this algorithm is demonstrated with extensive numerical experiments on the WDN demand maximization problem.

The rest of the paper is organized as follows. We begin by summarizing the literature on algorithms to solve water network optimization problems and feasible solution recovery for MILP-based approaches in Section \ref{sec:lit-rev}. Then, we present the statement of the WDN demand maximization problem and discuss the model's governing water flow equations in Section \ref{sec:statement-models}. Later, in Section \ref{sec:formulations}, we present the MINLP formulation to solve the water demand maximization problem and their MILP relaxations. Then, we will present the algorithm to recover feasible MINLP solutions from that of their MILP relaxations, which is the main contribution of the paper, in Section \ref{sec:algorithm}. Next, we discuss the performance of the algorithm in Section \ref{sec:results} through extensive numerical experiments performed on water networks studied in the literature. Finally, we summarize the work and discuss potential future directions in Section \ref{sec:summary}.

{
}

\section{Literature Review}
\label{sec:lit-rev}

The literature has studied several optimization problems in planning, design, and operations of water distribution networks. In addition to water demand maximization, these problems include optimal design or network expansion of the networks \cite{raghunathan2013global} and scheduling pumps to minimize the network operational costs \cite{tasseff2024polyhedral,bonvin2021pump}, among others. These problems are closely related to one other and most include aspects of the demand maximization problem. In their most general form, these problems are modeled as a MINLP.

Typically, generic MINLP solvers such as Juniper \cite{juniper}, BARON \cite{zhang2024solving}, and SCIP \cite{BestuzhevaEtal2021OO} fail to converge to an optimal solution or even a feasible solutions for these problems on medium to large-scale networks. To address this computational challenge, researchers developed several tailor-made algorithms \cite{menke2016exploring, raghunathan2013global,bonvin2021pump,costa2016branch,shi2016energy,naoum2015simulation,ghaddar2015lagrangian,bragalli2012optimal,nannicini2008local} over the past decade. These algorithms focus on computing global optimal solutions and utilize a range of techniques such as branching and bounding \cite{bonvin2021pump,nannicini2008local} to decomposition methods \cite{raghunathan2013global, ghaddar2015lagrangian} to simulation-based optimization \cite{naoum2015simulation} methods.

Recently, due to the superiority of the modern-day MILP solvers, MILP-based methods are gaining traction in solving MINLPs, especially those related to fluid-transport-infrastructure planning problems \cite{tasseff2024polyhedral, bonvin2021pump, hari2023relaxations}. These algorithms construct linear or piece-wise linear relaxations \cite{geissler2011mixed,tasseff2024polyhedral,bonvin2021pump} to the nonlinear constraints of the MINLP to form MILP relaxations.
Then, they embed the relaxations in either a branch-and-bound framework \cite{bonvin2021pump, nannicini2008local,sherali1997global} or a two-stage decomposition framework \cite{shi2016energy} to compute global optimal solutions.
In these frameworks, the relaxations are used to generate either bounds on the optimal objective values or candidate partial solutions and cuts to help the search.

While these algorithms focus on computing global optimal solutions, they fall short when a quick feasible solution is required.
Even upon early termination, they are not necessarily suited to provide a feasible solution or a good-quality solution as they get caught in the tradeoff between relaxation quality and their solve time.
Algorithms that use fixed partitions to generate relaxations \cite{bonvin2021pump} face the issue of either weak relaxations providing weak bounds and cuts (and therefore increasing the number of iterations required) or tight relaxations with each step of the iteration requiring long computation time  \cite{bragalli2012optimal}. Algorithms that allow partition refinement \cite{shi2016energy} do not focus on repairing infeasible candidate solutions and instead refine the partitions to generate new candidate solutions (with the goal of finding optimal solutions).

In this work, we focus on computing quick feasible solutions. This work differs from the aforementioned ones as we consider both relaxation (partition) refinement and also feasible solution recovery from infeasible solutions provided by a relaxations of \textit{arbitrary} partition size.
Consequently, the algorithm proposed here can also be incorporated into the aforementioned ones to accelerate the solution search.




The idea of computing feasible solutions with the aid of simpler versions of the problem is not new \cite{TELES2008376, TELES20093736, bonvin2021pump, nannicini2008local}.
For example, \cite{TELES2008376,TELES20093736} consider LP/MILP-based initialization of NLP solvers. The authors of \cite{fischetti2005feasibility, d2012storm} propose the idea of feasibility pumps to find integer feasible solutions for Mixed Integer Programs (MIPs) from continuous relaxations using techniques such as rounding.
The simulation-based branch and bound algorithm presented in \cite{bonvin2021pump} proposes a repair step for near-feasible NLP solutions; this step involves changing the time-discretization levels used to formulate the problem (thereby modifying the exact problem solved) to correct violations in the tank-levels of the infeasible solution.
The authors of \cite{ghaddar2015lagrangian} consider a Lagrangian decomposition method, where they apply a simulation-based limited discrepancy search to convert solutions of the Lagrangian relaxation of the problem into feasible solutions to the original problem.

While the limited discrepancy search of \cite{ghaddar2015lagrangian} has some similarities to the neighborhood search considered in our work, we focus on recovering solutions from MILP-relaxations instead of Lagrangian-relaxed solutions. 
Besides, in \cite{ghaddar2015lagrangian}, the authors progressively fix each binary variable, verify the feasibility of a few operational constraints, and perform a simulation using EPANET once all the binary variables are fixed. Contrastingly, in our work, we simultaneously search over all the integer variables by solving a restricted MINLP, as we will discuss in Section \ref{sec:algorithm}. 
Unlike the outer approximation used in \cite{shi2016energy} to create MILPs, here, we use outer and inner linear \textit{relaxations} to construct the MILPs.

More importantly, we combine sequences of converging piecewise linear relaxations with iterative relaxation refinement and feasible solution recovery from relaxations constructed from an \textit{arbitrary} number of partitions to develop a unique algorithm that quickly computes feasible MINLP solutions which progressively approach the optimal ones with time.



\section{Problem Statement \& Water Flow Models}
\label{sec:statement-models}
In this section, we formalize the water demand maximization problem for a WDN comprising reservoirs, pipes, pumps, tanks, and demand points. 
\subsection{Problem Statement}
Given a time series of maximum withdrawal flow rates required at each demand point within a specified time horizon, the problem objective is to determine the optimal operation of the network components, i.e., their on/off status and water flow rates, that maximize the total demand satisfied across all demand points and time points. The flow through each component is constrained by the physical laws governing the flow of water and operational limits.



    %




\subsection{Mathematical Models for Water Flow Through Network Components}
\label{subsec:water-flow-models}
\subsubsection{Junctions: Head Variables}
Water flow between two distinct locations in a network depends on the presence of a connecting component between the locations and is driven by the differences in total hydraulic heads across the components. The total hydraulic head of water at a given location is the sum of its pressure and elevation heads and is referred to as the \emph{head}. To model water flow across a network, we first define water head at connecting component intersections, referred to as junctions. We denote the set of all junctions in the network by $\mathcal{J}$. Then, for each junction, $j \in \mathcal J$, and every time point $t \in \mathcal{T}$,  water head is modeled with the variable $h_{j,t}$. The head at a junction is limited by factors such as the elevation of the junction and water pressure limits, as modeled  through the constraint set \eqref{junction:head_bounds},

\begin{align}
    \mathbf{\underline{h}_j} \leq h_{j,t} \leq \mathbf{\overline{h}_j}, \quad \forall j \in \mathcal{J}, t \in \mathcal{T},
    \label{junction:head_bounds}
\end{align}
where $\mathbf{\underline{h}_j}$ and $\mathbf{\overline{h}_j}$ denote the lower and upper bounds on the head at junction $j$, respectively.
Utilizing this fundamental variable, we next model the flow across several network connecting components.

\subsubsection{Pipes}
Pipes carry water between two junctions, with flow driven by the head difference between these junctions. 
We model pipes as edges and denote the set of all pipes by $\mathcal{A}_p$.
For every pipe, $a \in \mathcal{A}_p$, we designate an arbitrary `from' and `to' junctions and denote them by fr(a) and to(a) respectively, for ease of reference. 
We represent the net water flow rate from $fr(a)$ to $to(a)$ by the variable $q_{a,t}$. Note that the sign of $q_{a,t}$ is unrestricted since the junctions $fr(a)$ and $to(a)$ do not necessarily represent the true water flow direction. In fact, the flow direction may flip over the time horizon.

In hydraulic network optimization, two formulations are commonly used to model head loss along pipes \cite{tasseff2024polyhedral, bonvin2021pump}: the implicit-direction formulation, which expresses head loss as a function of 
$q_{a,t}|q_{a,t}|$, and the explicit-direction formulation, which introduces binary direction variables and nonnegative flow rates. In this work, we adopt the latter because it readily enables piecewise-linear relaxations and provides additional advantages to the algorithm developed, as discussed later.

Adopting the explicit-direction formulation, we model the water flow direction along a pipe $a$ using a binary decision variable $y_{a,t}$ which takes the value 1 if the flow is from $fr(a)$ to $to(a)$ and 0 otherwise.
Then, we represent the water flow rate along each pipe direction using two non-negative variables $q_{a,t}^+$ and $q_{a,t}^-$. One of $q_{a,t}^+$ and $q_{a,t}^-$ equals the volumetric flow rate depending on whether $y_{a,t}$ equals 1 or 0 respectively, and the other takes the value 0. 
We model the relationship between the flow direction, flow rates, and their limits using the set of constraints \eqref{pipe:flow_bounds-b} -- \eqref{pipe:flow_bounds-e}, 

\begin{align}
    \label{pipe:flow_bounds-b}
    q_{a,t} =  q_{a,t}^+ - q_{a,t}^-, \quad \forall a \in \mathcal{A}_p, t \in \mathcal{T},\\
    q_{a,t}^+ \leq \mathbf{\overline{q}_{a}^+} y_{a,t}, \quad \forall a \in \mathcal{A}_p, t \in \mathcal{T},\\
    q_{a,t}^+ \geq \mathbf{\underline{q}_{a}^+} y_{a,t}, \quad \forall a \in \mathcal{A}_p, t \in \mathcal{T},\\
    q_{a,t}^- \leq \mathbf{\overline{q}_{a}^-} (1 - y_{a,t}), \quad \forall a \in \mathcal{A}_p, t \in \mathcal{T},\\
    q_{a,t}^- \geq \mathbf{\underline{q}_{a}^-} (1-y_{a,t}), \quad \forall a \in \mathcal{A}_p, t \in \mathcal{T},
    \label{pipe:flow_bounds-e}
\end{align}
where $\mathbf{\overline{q}_{a}^+}$, $\mathbf{\overline{q}_{a}^-}$, denote the upper bounds and $ \mathbf{\underline{q}_{a}^+}$, and $ \mathbf{\underline{q}_{a}^-}$ denote the lower bounds on the flow rate along the indicated directions.

Next, recall that the head difference across pipes drives water flow through the pipes.
For each pipe $a$, we use the decision variables, $h_{a,t}^+ $ and $h_{a,t}^-$, to represent the non-negative head difference between the pipe endpoints in each direction (consistent with the above + and - directions). 
Then, the relationship between the head difference variables, the head variables across the pipe's `from' and `to' junctions, i.e., $h_{fr(a),t}$ and $h_{to(a),t}$ respectively, and the water flow direction, are modeled using constraints \eqref{pipe:head-b} -- \eqref{pipe:head-e}.

\begin{align}
    \label{pipe:head-b}
    \Delta h_{a,t}^+ \leq \mathbf{\overline{\Delta h}_{a}^+} y_{a,t}, \quad \forall a \in \mathcal{A}_p, t \in \mathcal{T},\\
    \Delta h_{a,t}^- \leq \mathbf{\overline{\Delta h}_{a}^-} (1 - y_{a,t}), \quad \forall a \in \mathcal{A}_p, t \in \mathcal{T},\\
    \Delta h_{a,t}^+ - \Delta h_{a,t}^- = h_{fr(a),t} - h_{to(a),t}, \quad \forall a \in \mathcal{A}_p, t \in \mathcal{T},
    \label{pipe:head-e}
\end{align}
where $\mathbf{\overline{\Delta h}_{a}^+}$ and $\mathbf{\overline{\Delta h}_{a}^-}$ denote the upper bounds on the head difference in each direction. When not given, these bounds are computed from the maximum and minimum heads achievable at each junction.

Given these pipe variables, the relationship between the flow rate along each direction of pipe $a$ and the head difference that drives flow is modeled with
\begin{align}
\label{pipe:head_loss-b}
    \Delta h_{a,t}^+ = \mathbf{r_a L_a} (q_{a,t}^+)^{1.852},\quad \forall a \in \mathcal{A}_p, t \in \mathcal{T},\\
    \Delta h_{a,t}^- = \mathbf{r_a L_a} (q_{a,t}^-)^{1.852}, \quad \forall a \in \mathcal{A}_p, t \in \mathcal{T},
    \label{pipe:head_loss-e}
\end{align}
where $L_a$ and $r_a$ denote the pipe's length and resistance per unit length, respectively.
Constraints \eqref{pipe:head_loss-b} and \eqref{pipe:head_loss-e} represent the well-known Hazen-Williams equation \cite{ormsbee2016darcy}, and the constant $r_a$ comprises all the length-independent constants that appear in the Hazen-Williams equation.


\subsection{Pumps}
Pumps are essential in boosting water head to compensate for the head loss along pipes and meet the delivery requirements at the demand nodes.
Here, we consider fixed-speed pumps that permit a \textit{unidirectional flow} and provide a head boost along the flow direction. 
We model these pumps as edges and represent the set of all pumps by $\mathcal{A}_{pu}$.

Pumps are controllable components that can be switched on or off to meet operational requirements. 
Generally, activating (switching on) a pump requires a certain minimum positive flow to avoid overheating. If the flow is zero, pumps are inactive (switched off).
Here, we model the activation (on/off) status of pump $a \in \mathcal{A}_{pu}$ using the binary variable $z_{a,t}$, which takes the value 1 if the pump is active and 0 otherwise. 

Like pipes, we denote the junctions at which water enters and leaves pump $a$ by $fr(a)$ and $to(a)$, respectively. Since the flow is unidirectional, these junctions also represent the flow direction. Next, the water flow rate along pump $a$ from $fr(a)$ to $to(a)$ is represented by the non-negative variable $q_{a,t}$.
Then, we model the relationship between the pump's activation status and water flow rate with:
    \begin{align}
        \label{pump:flow_b}
        q_{a,t} \geq \mathbf{\underline{q}_{a}} z_{a,t}, \quad \forall a \in \mathcal{A}_{pu}, t \in \mathcal{T},\\
        q_{a,t} \leq \mathbf{\overline{q}_{a}} z_{a,t}, \quad \forall a \in \mathcal{A}_{pu}, t \in \mathcal{T},
        \label{pump_flow_e}
    \end{align}
     where $\mathbf{\underline{q}_{a}}$ and $\mathbf{\overline{q}_{a}}$ represent the positive lower and upper bounds on the flow rate when the pump is active.

Pumps provide a head boost only when active. Otherwise, the heads at the junctions connected by the pump are decoupled and there is no flow. The head gain offered by pump $a$ is modeled with the non-negative variable $g_{a,t}$. Then, we express the relationship between a pump's activation status, the head gain the pump provides, and the heads across the junctions connected by the pump using the following constraints.

    \begin{multline}
        h_{to(a),t} - h_{fr(a),t} \leq g_{a,t} + \mathbf{\overline{\Delta h}_{a}^-} (1 - z_{a,t}),\\ \quad \forall a \in \mathcal{A}_{pu}, t \in \mathcal{T},
        \label{pump:head-b}
    \end{multline}
    \begin{multline}
        h_{to(a),t} - h_{fr(a),t} \geq g_{a,t} +  \mathbf{\underline{\Delta h}_{a}^-} (1 - z_{a,t}),\\ \quad \forall a \in \mathcal{A}_{pu}, t \in \mathcal{T},
        \label{pump:head-e}
    \end{multline}
    where $\mathbf{\overline{\Delta h}_{a}^-}$ and $\mathbf{\underline{\Delta h}_{a}^-}$ represent the upper and lower bounds on the permitted head gain as defined by the head limits at the junctions.

The value of the head gain offered by the pump is typically determined from head curves provided by the pump's manufacturer. To simplify the incorporation of these curves into water distribution planning, researchers have developed approximations \cite{ulanicki2008modeling, coulbeck1991ginas, rossman2000epanet}. Here, we use one such approximation \cite{ulanicki2008modeling}, where the head gain is a strictly concave quadratic function of the water flow rate and is non-zero only when the pump is active. We express this mathematically as shown in the constraint set \eqref{pump:head_gain}.

    \begin{align}
        \boldsymbol{\alpha_a} (q_{a,t})^2 + \boldsymbol{\beta_a} (q_{a,t}) + \boldsymbol{\gamma_a} z_{a,t} = g_{a,t}, \forall a \in \mathcal{A}_{pu}, t \in \mathcal{T}
        \label{pump:head_gain}
    \end{align}
    where the coefficients $\boldsymbol{\alpha_a}$, $\boldsymbol{\beta_a} $, and $\boldsymbol{\gamma_a}$ obey the conditions $\boldsymbol{\alpha_a} < 0$ and $\boldsymbol{\gamma_a} > 0$, and are calculated from the pump curve.  

Finally, pumps consume energy to provide head gain. The power consumption of pumps is typically a nonlinear function of water flow rate, efficiency, and head gain \cite{ulanicki2008modeling}. However, fixed-speed pumps are often approximated with an affine function. Specifically, the power consumption of pump $a \in \mathcal{A}_{pu}$ is modeled as an affine function of its water flow rate and activation status, as shown in \eqref{pump:power_consumption}.
\begin{align}
    P_{a,t} = \boldsymbol{\omega_a}  q_{a,t} + \boldsymbol{\mu_a} z_{a,t}, \quad \forall a \in \mathcal{A}_{pu}, t \in \mathcal{T},
    \label{pump:power_consumption}
\end{align}
{where $\boldsymbol{\omega_a} $ and $\boldsymbol{\mu_a}$ are determined from pump power curves.}

\subsection{Tanks}              
Tanks help in storing and discharging water to meet demand when required (usually at peak hours). Here, we denote the set of all tanks in the network by $\mathcal{TK}$ and model every tank $i \in \mathcal{TK}$ as a cylindrical storage element with a constant cross-sectional area, $\mathbf{A_{i}}$.
We assume that the tanks are vented to the atmosphere and there is no pressure head.

Every tank, $i \in \mathcal{TK}$, is connected to a junction, denoted by $j_i \in \mathcal{J}$. We assume that the elevation of the tank bottom, $\mathbf{b_{j_i}}$, is no more than the minimum head at the junction, $\mathbf{\underline{h_{j_i}}}$.
Then, the water volume in the tank, denoted by the variable, $V_{i,t}$, is related to the head through the constraints:

\begin{align}
    V_{i,t} = \mathbf{A_i} (h_{j_i,t} - \mathbf{b_{i}}), \quad \forall i \in \mathcal{TK}, t \in \mathcal{T}.
    \label{tank:volume-head}
\end{align}

At the beginning of the time horizon, the volume of water present in the tank is denoted by $V_{i,0}$. For the remaining horizon, the water volume must remain within specified upper and lower bounds, $\mathbf{\overline{V_i}}$ and $\mathbf{\underline{V_i}}$ respectively, which are determined by the tank capacity and other operational limits. 
These conditions and volume bounds are enforced using the following constraints.

    \begin{align}
        \label{tank:initial_volume}
        V_{i,T_1} = \mathbf{V_{i,0}}, \quad \forall i \in \mathcal{TK}\\
        \mathbf{\underline{V_i}}\leq V_{i,t} \leq \mathbf{\overline{V_i}}, \quad \forall i \in \mathcal{TK}, t \in \mathcal{T}.
        \label{tank_volume_bounds}
    \end{align}


After each time step of duration $\mathbf{\Delta t}$, the tank water volume changes based on the net withdrawal or injection. Let the variable $q_{i,t}$ denote the net withdrawal rate from tank $i$ at time $t$. Then, the volume update is enforced by constraint set \eqref{tank:volume_update}:
\begin{align}
    V_{i,\,t+\mathbf{\Delta t}} 
    &= V_{i,t} - q_{i,t}\,\mathbf{\Delta t}, 
    \quad \forall i \in \mathcal{TK}, \; t \in \mathcal{T}\setminus\{\mathbf{T_f}\},
    \label{tank:volume_update}
\end{align}
where $\mathbf{T_f}$ is the final time point in the horizon.

\subsection{Reservoirs}
Similar to \cite{tasseff2024polyhedral}, we assume reservoirs are infinite sources of supply with a constant head (zero pressure head and constant elevation) that inject water into the distribution network at any flow rate. 
We model reservoirs as elements connected directly to a junction and denote the set of all reservoirs by $\mathcal{R}$.
We denote the rate of injection from reservoir $i \in \mathcal{R}$ with the non-negative variable $q_{i,t}$. Non-negativity indicates that the flow is unidirectional, and we express this as \eqref{reservoir:non-neg}.

\begin{align}
    q_{i,t} \geq 0,  \quad \forall i \in \mathcal{R}, t \in \mathcal{T}
    \label{reservoir:non-neg}
\end{align}

\subsection{Demand Points}
Unlike reservoirs, demand points are points where water is withdrawn from the network. The rate of withdrawal is capped by the specified maximum demand required by the consumers at these points. Nevertheless, similar to reservoirs, we model them as elements directly connected to junctions at which the withdrawal takes place.
We denote the set of all demand points by $\mathcal{D}$ and the variable flow rate of withdrawal at demand point $i \in \mathcal{D}$ by $q_{i,t}$. Then, the withdrawal cap is modeled by constraints \eqref{demand:withdrawal-cap},

 \begin{align}
     q_{i,t} \leq \mathbf{d_{i,t}}, \quad \forall i \in \mathcal{D}, t \in \mathcal{T},
     \label{demand:withdrawal-cap}
 \end{align}   
 where $\mathbf{d_{i,t}}$ is the maximum withdrawal rate required at the demand point $i$ and time point $t$.


\subsection{Junctions: Flow Conservation}
Everywhere within the network, water flow must satisfy the law of mass conservation. To enforce this, at every junction, we require the water inflow rate to be equal to water outflow rate.
Given a junction $j \in \mathcal{J}$, let  $\mathcal{C}_j^+$ and $\mathcal{C}_j^-$  denote the sets of all components (pipes, pumps, tanks, reservoirs, and demand points) with flows into and out of the junction, respectively. Then, mass flow balance is enforced by constraints  \eqref{junction:mass_conservation},

    \begin{align}
        \sum_{a \in \mathcal{C}_j^+}  q_{a,t} = \sum_{b \in \mathcal{C}_j^-} q_{b,t}, \quad \forall j \in \mathcal{J}, t \in \mathcal{T}.
        \label{junction:mass_conservation}
    \end{align}



\subsection{Objective}
Finally, while satisfying the above-discussed physics, the problem objective is to maximize the sum of water flow rates supplied to all the demand points across all the time points. Conforming to the notation, we can express the objective function that must be maximized as

\begin{equation}
    \sum_{d \in \mathcal{D},t \in \mathcal{T}} q_{d,t}
    \label{minlp:obj}
\end{equation}


\section{Mathematical Formulations}
\label{sec:formulations}
\subsection{Problem Formulation: Mixed Integer Nonlinear Program (MINLP)}
\label{subsec:MINLP}
Here, the above constraints are combined to define the demand maximization (DM) problem as
\begin{equation} \label{form:minlp}
    \begin{aligned}
     (\mathcal N_1)
        & \text{ Max} & & \begin{tabular}{@{}c@{}} \textnormal{Total water demand satisfied: } \end{tabular}  \eqref{minlp:obj},\\
        & \text{s.t.} & & \textnormal{Pipe constraints: } \eqref{pipe:flow_bounds-b} - \eqref{pipe:head_loss-e}\\
        & & & \textnormal{Pump constraints: } \eqref{pump:flow_b} - \eqref{pump:head_gain},\\
        & & & \textnormal{Tank constraints: } \eqref{tank:volume-head} - \eqref{tank:volume_update},\\
                & & & \textnormal{Reservoir: } \eqref{reservoir:non-neg},\\
        & & & \textnormal{Demand constraints: } \eqref{demand:withdrawal-cap}\\
        & & & \textnormal{Junction constraints: } \eqref{junction:head_bounds}, \eqref{junction:mass_conservation}\\
    \end{aligned}
\end{equation}

\subsection{Problem Relaxation: Mixed Integer Linear Program (MILP) Relaxation}
\label{subsec:milp-relaxation}
In the formulation $\mathcal N_1$, nonlinearities in the continuous space arise from the head-loss and head-gain constraints~\eqref{pipe:head_loss-b}–\eqref{pipe:head_loss-e} and~\eqref{pump:head_gain}. 
Let $\Omega$ denote the domain of the pipe-flow and pump-flow variables. 
To relax these nonlinearities, we partition $\Omega$ into intervals and construct piecewise-linear inner and outer envelopes for  constraints~\eqref{pipe:head_loss-b}, \eqref{pipe:head_loss-e} and~\eqref{pump:head_gain}, following the procedure described in~\cite{tasseff2024polyhedral}. 

\begin{definition}[Partition]
Given a variable domain $\Omega$, its \emph{partition} $\mathcal P$ is defined as a finite collection of intervals
\[
\mathcal P = \{\Delta_1, \Delta_2, \ldots, \Delta_m\},
\qquad
\Omega = \bigcup_{i=1}^m \Delta_i,
\]
where adjacent intervals may share endpoints but have no overlapping interiors. 
Each $\Delta_i$ is referred to as a \emph{partitioned interval}.
\end{definition}

For a given partition $\mathcal P$, we use the convex-combination method \cite{tasseff2024polyhedral} to formulate the piecewise-linear envelopes and construct an MILP relaxation of $\mathcal N_1$\footnote{The detailed construction of these relaxations is well established in the literature; see~\cite{tasseff2024polyhedral, sundar2022sequence} for the complete procedure.}.
Let $\mathcal{HW}_r(\mathcal P)$ and $\mathcal{HG}_r(\mathcal P)$ denote the constraint sets representing the piecewise relaxations of the pipe head-loss constraints~\eqref{pipe:head_loss-b}–\eqref{pipe:head_loss-e} and the pump head-gain constraints~\eqref{pump:head-b}–\eqref{pump:head-e}, respectively. 
Then, the resulting MILP relaxation of $\mathcal N_1$ is denoted $\mathcal L_1(\mathcal P)$ and is formulated as:

\begin{equation} \label{form:milp}
    \begin{aligned}
     (\mathcal L_1 (\mathcal{P}))
        & \text{ Max} & & \begin{tabular}{@{}c@{}} \textnormal{Total water demand satisfied: } \end{tabular}  \eqref{minlp:obj},\\
        & \text{ s.t.} & & \textnormal{Pipe constraints: } \eqref{pipe:flow_bounds-b} - \eqref{pipe:head-e}, \mathcal{HW}_r(\mathcal{P})\\
        & & & \textnormal{Pump constraints: } \eqref{pump:flow_b} - \eqref{pump:head-e}, \mathcal{HG}_r(\mathcal{P})\\ 
        & & & \textnormal{Tank constraints: } \eqref{tank:volume-head} - \eqref{tank:volume_update},\\
                & & & \textnormal{Reservoir: } \eqref{reservoir:non-neg},\\
        & & & \textnormal{Demand constraints: } \eqref{demand:withdrawal-cap}\\
        & & & \textnormal{Junction constraints: } \eqref{junction:head_bounds}, \eqref{junction:mass_conservation}\\
    \end{aligned}
\end{equation}

$\mathcal L_1 (\mathcal{P})$ may sometimes exhibit solution degeneracy where there are multiple solutions with the same objective value, which can lead to undesirable algorithmic trends and arbitrary solution selection. To break solution equivalence systematically in such cases we introduce a second MILP,  $\mathcal L_2 (\mathcal{P})$. In this MILP, we minimize the sum of pump operational costs and the reservoir supply rate over the time horizon by fixing the demand variables to the solutions of $\mathcal L_1 (\mathcal{P})$.

Suppose $\{\mathbf{q^*_{i,t}}: \forall i \in \mathcal{D}, t \in \mathcal{T}\}$ is the set of demand solutions obtained by solving $\mathcal L_1 (\mathcal{P})$. Then, in $\mathcal L_2 (\mathcal{P})$, we fix these demand values by replacing the demand constraints \eqref{demand:withdrawal-cap} with the following constraints

\begin{align}
    q_{i,t} = \mathbf{q^*_{i,t}}, \quad \forall i \in \mathcal{D}, t \in \mathcal{T}.
    \label{fixing_demand}
\end{align}

The objective function is then defined as \eqref{milp2:obj},
\begin{equation}
    \sum_{t \in \mathcal{T}} (\sum_{i \in \mathcal{R}} q_{i,t} + \sum_{a \in \mathcal{A}_{pu}} (\boldsymbol{\$_{a,t} \Delta t} ) P_{a,t}),
    \label{milp2:obj}
\end{equation}
where $\boldsymbol{\$_{a,t}}$ represents the energy price for pump $a$ at time point $t$.
Then, the tie-breaking MILP is expressed as:
\begin{equation} \label{form:milp-2}
    \begin{aligned}
     (\mathcal L_2 (\mathcal{P}))
        & \text{ Min} & & \begin{tabular}{@{}c@{}} \textnormal{Pump costs and reservoir supply: } \end{tabular}  \eqref{milp2:obj},\\
        & \text{ s.t.} & & \textnormal{Pipe constraints: } \eqref{pipe:flow_bounds-b} - \eqref{pipe:head-e}, \mathcal{HW}_r(\mathscr{P})\\
        & & & \textnormal{Pump constraints: } \eqref{pump:flow_b} - \eqref{pump:head-e}, \mathcal{HG}_r(\mathscr{P})\\ 
        & & & \textnormal{Tank constraints: } \eqref{tank:volume-head} - \eqref{tank:volume_update},\\
                & & & \textnormal{Reservoir: } \eqref{reservoir:non-neg},\\
        & & & \textnormal{Fixed Demand constraints: } \eqref{fixing_demand}\\
        & & & \textnormal{Junction constraints: } \eqref{junction:head_bounds}, \eqref{junction:mass_conservation}
    \end{aligned}
\end{equation}


\color{black}

\section{Algorithm}
\label{sec:algorithm}

In this section, we first present the baseline recovery algorithm, 
a standard approach for computing feasible solutions to $\mathcal N_1$ using solutions of its MILP relaxations such as $\mathcal{L}_1(.)$.
We then discuss its shortcomings and progressively enhance the algorithm, leading to a robust framework for quickly recovering high-quality feasible MINLP solutions. For ease of description, let $\mathcal X = (x_1, \ldots, x_{\mathbf{n_I}})$ and $\mathcal Y = (y_1, \ldots, y_{\mathbf{n_C}})$ denote the ordered tuples of the integer and continuous variables of $\mathcal{N}_1$, respectively, where $\mathbf{n_I}$ and $\mathbf{n_C}$ are the respective tuple sizes.
Let $c(\mathcal S)$ denote the objective value of a feasible solution $\mathcal S = (\mathcal X, \mathcal Y)$ of $\mathcal N_1$, and $c_{\mathcal L_1(\mathcal P_k)}^\star$ and $c_{\mathcal L_2(\mathcal P_k)}^\star$ denote the optimal objective values of the MILP relaxations $\mathcal L_1(\mathcal P_k)$ and $\mathcal L_2(\mathcal P_k)$, respectively.

\subsection{Baseline Recovery Algorithm}
Suppose that a partition, $\mathcal P$, of $\Omega$ is given to construct and solve $\mathcal{L}_1(\mathcal P)$, an MILP relaxation to $\mathcal N_1$. 
Let $\Tilde{\mathcal X} = (\Tilde{x}_1, \ldots, \Tilde{x}_{\mathbf{n_I}})$ denote the components of the computed optimal solution of $\mathcal{L}_1(\mathcal P)$ corresponding to the variable tuple $\mathcal X$. Then, the baseline algorithm for recovering a feasible solution to $\mathcal N_1$ from the candidate solution $\Tilde{\mathcal X}$ is presented in Algorithm \ref{alg:baseline}.

\begin{algorithm}[H]
\caption{Baseline Recovery Algorithm}
\label{alg:baseline}
\begin{algorithmic}[1]
\Require $\mathcal{N}_1$, relaxation $\mathcal L_1(\mathcal P)$ \& candidate $\mathcal{\Tilde{X}}$.
\Ensure Feasible solution to $\mathcal{N}_1$
\State Form the NLP, $\mathcal{N}_1^{\tilde{\mathcal{X}}}$, by fixing $\mathcal{X}\gets\tilde{\mathcal{X}}$ in $\mathcal{N}_1$.
\State Solve $\mathcal{N}_1^{\tilde{\mathcal{X}}}$
\If{$\mathcal{N}_1^{\tilde{\mathcal{X}}}$ is feasible}
  \State obtain its optimal solution, $\tilde{\mathcal{Y}}^{\,\tilde{\mathcal{X}}}$ and \Return $(\tilde{\mathcal{X}},\,\tilde{\mathcal{Y}}^{\,\tilde{\mathcal{X}}})$
\Else
  \State \Return \textsc{Fail}
\EndIf
\end{algorithmic}
\end{algorithm}

The algorithm fixes the integer variables of $\mathcal{N}_1$ to their corresponding values in $\Tilde{\mathcal{X}}$, yielding the nonlinear program (NLP) $\mathcal{N}_1^{\Tilde{\mathcal{X}}}$. If $\mathcal{N}_1^{\mathcal{\Tilde{X}}}$ is feasible, solving it to optimum yields a continuous solution $\mathcal{\Tilde{Y}}^{\mathcal{\Tilde{X}}}$. Then, $(\mathcal{\Tilde{X}}, \mathcal{\Tilde{Y}}^{\mathcal{\Tilde{X}}})$ is the desired feasible solution to $\mathcal{N}_1$. However, if $\mathcal{N}_1^{\mathcal{\Tilde{X}}}$ is infeasible, the algorithm fails to produce a feasible solution to  $\mathcal{N}_1$.

The baseline algorithm often fails when initiated with a coarse partition $\mathcal P$, since the candidate $\Tilde{\mathcal X}$ obtained from $\mathcal L_1(\mathcal P)$ is typically infeasible for $\mathcal N_1$. Conversely, overly fine partitions make $\mathcal L_1(\mathcal P)$ intractable. To recover feasible solutions for any given partition $\mathcal P$, we extend the baseline method by incorporating a neighborhood search around $\Tilde{\mathcal X}$ and progressively enlarging the neighborhood until a feasible solution is found.

\subsection{Improvement to Baseline with Neighborhood Search}
We define the neighborhood size by the Hamming distance, $h$, i.e., the number of component-wise deviations permitted between the integer variable tuple $\mathcal X$ and the given candidate $\tilde{\mathcal X}$. This distance is enforced by adding constraint \eqref{eq:nbhd-hops} to $\mathcal N_1$, yielding a restricted MINLP, denoted $\mathcal{N}_2 (\mathcal{\Tilde{X}}, h)$.

\begin{equation}
\sum_{\{i \in \mathcal I: \Tilde{\boldsymbol{x}}_i = 0\}} x_i
\;+\;
\sum_{\{i \in \mathcal I: \Tilde{\boldsymbol{x}}_i = 1\}} (1 - x_i)
\;=\; h,
\label{eq:nbhd-hops}
\end{equation}
where $\mathcal I = \{1, \ldots, \mathbf{n_I}\}$ denotes the set of component indices of $\mathcal{X}$.
Then, the recovery algorithm based on iteratively enlarging neighborhood search is presented in Algorithm \ref{alg:baseline+hopping}.

\begin{algorithm}[H]
\caption{Recovery Algorithm with Neighborhood Search}
\label{alg:baseline+hopping}
\begin{algorithmic}[1]
\Require $\mathcal N_1$, $\mathcal L_1(\mathcal P)$ \& candidate $\mathcal{\Tilde{X}}$, time limit $\mathbf{T_{max}}$.
\Ensure Feasible solution to $\mathcal N_1$
\For{$h \gets 0$ \textbf{to} $\mathbf{n_I}$}
  \State Build $\mathcal N_2(\Tilde{\mathcal X}, h)$ by adding constraint \eqref{eq:nbhd-hops} to $\mathcal N_1$ \label{step:build n2}
  \State Solve $\mathcal N_2(\Tilde{\mathcal X}, h)$ with time limit $\mathbf{T_{max}}$
  \If{$\mathcal N_2(\Tilde{\mathcal X}, h)$ is feasible}
    \State Set $\bar h \gets h$
    \State \Return $(\mathcal X^{\bar h}, \mathcal Y^{\bar h})$, the best feasible solution of $\mathcal N_2(\Tilde{\mathcal X}, \bar h)$
    \Else{ Set $h \gets h + 1$; return to step \ref{step:build n2}}  
  \EndIf
\EndFor
\end{algorithmic}
\end{algorithm}

The algorithm proceeds by solving $\mathcal{N}_2 (\mathcal{\Tilde{X}}, h)$ iteratively for increasing values of $h$, starting from $h = 0$ and increasing it by $1$ until $\mathcal{N}_2 (\mathcal{\Tilde{X}}, h)$ is feasible. Let $\bar h$ denote the smallest value such that $\mathcal{N}_2 (\mathcal{\Tilde{X}}, \bar h)$ is feasible. Within the time limit $\mathbf{T_{max}}$, let $(\mathcal X^{\bar h},\mathcal Y^{\bar h})$ denote the best feasible solution returned by the solver for $\mathcal N_2(\tilde{\mathcal X}, \bar h)$. Then $(\mathcal X^{\bar h},\mathcal Y^{\bar h})$ is the desired feasible solution to $\mathcal N_1$.

Note that the algorithm generalizes the baseline method, since they are equivalent when $h = 0$. More importantly, unlike the baseline, it is guaranteed to return a feasible solution for any given candidate $\mathcal{\Tilde{X}}$ or partition $\mathcal P$, as $\mathcal{N}_2(\cdot)$ reduces to $\mathcal N_1$ when $h = |\mathcal I|$.
Nevertheless, $\mathcal{N}_2(\cdot)$ may be computationally challenging, particularly when the search involves a large number of integer variables.
To enhance the computational performance of the algorithm, we restrict the neighborhood search to a subset of integer variables.

\subsubsection{Restriction to a subset}
Let $\mathcal I_s \subset \mathcal I$ denote the subset of indices of integer variables selected for the neighborhood search. 
The restricted version of the algorithm differs from Algorithm \ref{alg:baseline+hopping} in the construction of $\mathcal N_2(\Tilde{\mathcal X}, h)$. 
Specifically, in Step \ref{step:build n2-restricted-nbhd}, variables not selected for neighborhood search are fixed to their candidate values in $\Tilde{\mathcal X}$, and the Hamming-distance constraint is enforced only on $\mathcal I_s$ by replacing constraint  \eqref{eq:nbhd-hops} with constraint \eqref{eq:nbhd-hops-restricted}.
\begin{equation}
\sum_{\{i \in \mathcal I_s: \Tilde{\boldsymbol{x}}_i = 0\}} x_i
\;+\;
\sum_{\{i \in \mathcal I_s: \Tilde{\boldsymbol{x}}_i = 1\}} (1 - x_i)
\;=\; h.
\label{eq:nbhd-hops-restricted}
\end{equation}
The resulting algorithm is given in Algorithm \ref{alg:restricted-neighborhood}.


\begin{algorithm}[H]
\caption{Recovery Algorithm with Neighborhood Search Restricted to a Subset of Integers}
\label{alg:restricted-neighborhood}
\begin{algorithmic}[1]
\Require $\mathcal N_1$, $\mathcal L_1(\mathcal P)$ \& candidate $\Tilde{\mathcal X}$, restricted index set $\mathcal I_s$, time limit $\mathbf{T_{max}}$.
\Ensure Feasible solution to $\mathcal N_1$
\For{$h \gets 0$ \textbf{to} $|\mathcal I_s|$}
  \State Build $\mathcal N_2(\tilde{\mathcal X}, h)$ by fixing $x_i \gets\tilde{\boldsymbol{x}_i}, \forall i \in \mathcal I \setminus \mathcal I_s$  in $\mathcal{N}_1$ and 
  adding constraint \eqref{eq:nbhd-hops-restricted} \label{step:build n2-restricted-nbhd}.
  \State Solve $\mathcal N_2(\tilde{\mathcal X}, h)$ with time limit $\mathbf{T_{max}}$.
  \If{$\mathcal N_2(\tilde{\mathcal X}, h)$ is feasible}
    \State Set $\bar h \gets h$
    \State \Return $(\mathcal X^{\bar h},\mathcal Y^{\bar h})$, the best feasible solution of $\mathcal N_2(\tilde{\mathcal X}, \bar h)$
    \Else{ Set $h \gets h + 1$; return to step \ref{step:build n2}} 
  \EndIf
\EndFor
\end{algorithmic}
\end{algorithm}

The specification of $\mathcal I_s$ governs the tradeoff between computational effort and solution quality, and identifying an effective subset typically requires either domain knowledge or trial-and-error. We discuss the choice of $\mathcal I_s$ for the water demand maximization problem considered in this paper in Section \ref{sec:results}. Notably, Algorithm \ref{alg:restricted-neighborhood} reduces to Algorithm \ref{alg:baseline+hopping} when $\mathcal I_s = \mathcal I$ and to Algorithm \ref{alg:baseline} when $\mathcal I_s$ is empty. Therefore, we focus the subsequent discussion on Algorithm \ref{alg:restricted-neighborhood}

The solution quality of all the aforementioned algorithms depends critically on the partition $\mathcal P$, since it determines the quality of the input candidate $\Tilde{\mathcal X}$. Finer partitions typically yield stronger candidates that facilitate solution recovery but require solving large, computationally challenging MILPs. In contrast, coarse partitions generate weaker candidates that can be obtained quickly and subsequently improved through the proposed neighborhood search. Besides, the level of partition granularity sufficient to produce strong candidates is unknown a priori.
To navigate the tradeoffs between candidate quality, generation effort, and the potential for improvement via neighborhood search, we embed Algorithm \ref{alg:restricted-neighborhood} within an iterative partition refinement algorithm.

\subsection{Integration With Partition Refinement}
This approach leverages a pool of candidates generated across various partitions and initiates solution recovery for each candidate by applying the neighborhood search as necessary.
To formalize this procedure, we define the notion of partition refinement.


\begin{definition}[Partition Refinement]
Given two partitions $\mathcal P_i$ and $\mathcal P_j$ of the same domain $\Omega$,  
$\mathcal P_j$ is said to \emph{refine} $\mathcal P_i$ if for every partitioned interval $\Delta_j \in \mathcal P_j$,  
there exists an interval $\Delta_i \in \mathcal P_i$ such that $\Delta_j \subseteq \Delta_i$.
\end{definition}
The initial partition $\mathcal P_0$ may be chosen as the entire domain $\Omega$ or any of its refinements.
While the algorithm is valid for any refinement, in this work we adopt a simple bisection rule, where at each refinement step every interval in $\mathcal P_k$ is split into two equal halves, yielding $\mathcal P_{k+1}$.

The resulting iterative partition refinement algorithm is presented in Algorithm \ref{alg:partition-refinement}.

\begin{algorithm}[H]
\caption{Iterative partition refinement with neighborhood search.}
\label{alg:partition-refinement}
\begin{algorithmic}[1]
\Require $\mathcal N_1$, $\mathcal L_1(\cdot)$ \& candidate $\Tilde{\mathcal X}$, $\mathcal I_s$, time limit $\mathbf{T_{max}}$, maximum refinements $\mathbf{K_{\max}}$, tolerance $\varepsilon_{\text{opt}}$
\Ensure Feasible solution to $\mathcal N_1$ 
\State Initialize $\mathcal P_0$ \Comment{e.g., $\mathcal P_0 \gets \Omega$}
\State Initialize $k \gets 1$, best incumbent $\mathcal S^\star \gets \varnothing$, and its cost $c^\star \gets -\infty$ 

\While{$k \le K_{\max}$ \textbf{and} the incumbent $\mathcal S^\star$ is not certified optimal for $\mathcal N_1$}
    \State $\mathcal P_{k} \gets \text{Refine}(\mathcal P_{k-1})$ \Comment{e.g., bisect each interval}
    \State Solve $\mathcal L_1(\mathcal P_k)$ within $T_{\max}$; if an optimal solution is found, let $\Tilde{\mathcal X}_k$ be the projected candidate, $c_{\mathcal L_1(\mathcal P_k)}^\star$ be its objective value, and $t_k^{L_1}$ be the solve time.
    \State Apply Algorithm~\ref{alg:restricted-neighborhood} with $(\mathcal N_1,\mathcal L_1(\mathcal P_k), \Tilde{\mathcal X_k}, \mathcal I_s, \mathbf{T_{max}})$ to attempt recovery
    \If{Algorithm \ref{alg:restricted-neighborhood} returns a feasible solution $\mathcal S$ with objective $c(\mathcal S) > c^\star$}
        \State $\mathcal S^\star \gets \mathcal S$; \quad $c^\star \gets c(\mathcal S)$
        \If{$|c^\star - c_{\mathcal L_1(\mathcal P_k)}^\star| \le \varepsilon_{\text{opt}}$}
            \State $\mathcal S^\star$ is \text{certified optimal for} $\mathcal N_1$; 
        \EndIf
    \EndIf
    \State $k \gets k+1$
\EndWhile
\State \Return $\mathcal S^\star$
\end{algorithmic}
\end{algorithm}

The algorithm involves solving the relaxation $\mathcal L_1(\mathcal P_k)$ for a given partition to generate an integer candidate,  attempting recovery from the candidate, and refining the partition.
Iteration continues until termination by optimality—when the cost of the best incumbent solution is within a specified tolerance of the best upper bound $c_{\mathcal L_1(\mathcal P_k)}^\star$—or until the time limit $\mathbf{T_{\max}}$ or the refinement iteration limit $\mathbf{K_{\max}}$ is reached.

\subsection{Complete Algorithm with Tie-Breaker}
Generally, the objective value of the recovered solution is expected to improve or at least remain unchanged as the partition is refined. 
However, this trend is not always observed, since $\mathcal L_1(\cdot)$ may admit multiple optimal solutions that yield different integer candidates. 
Although these candidates share the same relaxation objective, they can lead to different recovered objective values across partitions, resulting in inconsistencies. 
To enforce consistency and systematically break ties among equivalent candidates, we employ the tie-breaking MILP $\mathcal L_2(\cdot)$.

Note that solving $\mathcal L_2(\cdot)$ can sometimes be more challenging than solving $\mathcal L_1(\cdot)$. 
Nevertheless, by construction, any feasible solution to $\mathcal L_2(\cdot)$ is also optimal for $\mathcal L_1(\cdot)$ and therefore valid for recovery. 
To maintain tractability, a time limit is imposed when solving $\mathcal L_2(\cdot)$, and the best available candidate within this limit is retained for recovery.

$\mathcal L_2(\cdot)$ not only resolves ties but also provides an additional candidate for recovery. To exploit this, 
To exploit the consistency and the availability of an additional candidate, we augment the previous algorithm with a recovery step from the $\mathcal L_2(\cdot)$ candidate.
The complete iterative refinement procedure is summarized in Algorithm \ref{alg:partition-refinement-tiebreak}.

\begin{algorithm}[ht]
\caption{Iterative partition refinement with neighborhood search and embedded tie-breaking recovery}
\label{alg:partition-refinement-tiebreak}
\begin{algorithmic}[1]
\Require $\mathcal N_1$, $\mathcal L_1(\cdot)$, $\mathcal L_2(\cdot)$, $\mathcal I_s$, time limit $\mathbf{T_{max}}$, maximum refinements $\mathbf{K_{\max}}$, tolerance $\varepsilon_{\text{opt}}$
\Ensure Feasible solution to $\mathcal N_1$
\State Initialize $\mathcal P_0$ \Comment{e.g., $\mathcal P_0 \gets \Omega$}
\State Initialize $k \gets 1$, best incumbent $\mathcal S^\star \gets \varnothing$, and its cost $c^\star \gets -\infty$ 
\While{$k \le K_{\max}$ \textbf{and} the incumbent $\mathcal S^\star$ is not certified optimal for $\mathcal N_1$}
\State $\mathcal P_k \gets \text{Refine}(\mathcal P_{k-1})$ \Comment{e.g., bisect each interval}
    \State \textbf{($\mathcal L_1$ step)} Solve $\mathcal L_1(\mathcal P_k)$ within $\mathbf{T_{max}}$; 
     if optimal, let $\Tilde{\mathcal X}^{(1)}_k$ be the projected candidate, 
           $\mathbf{\tilde q}^k$ the demand values, 
           $c_{\mathcal L_1(\mathcal P_k)}^\star$ its objective value, 
           and $t_k^{L_1}$ the solve time
    \label{a5:step-l1}
    \State Set the candidate generation time, $t_k^{LP} \gets t_k^{L_1}$
    \label{a5:step-l1-time}
    \State \textbf{(Recovery from $\mathcal L_1$)} Apply Algorithm~\ref{alg:restricted-neighborhood} with $(\mathcal N_1,\mathcal L_1(\mathcal P_k),\Tilde{\mathcal X}^{(1)}_k, \mathcal I_s,\mathbf{T_{max}})$; record the recovery time as $t_k^{rec1}$
    \label{a5:l1-rec-begin}
    \If{ Algorithm \ref{alg:restricted-neighborhood} returns a feasible solution $\mathcal S_1$}
        \State Set the solution recovery time, $t_k^{rec} \gets t_k^{rec1}$
         \If{ objective $c(\mathcal S_1)>c^\star$}
            \State $\mathcal S^\star \gets \mathcal S_1$;\quad $c^\star \gets c(\mathcal S_1)$
        \EndIf
    \EndIf
    \label{a5:l1-rec-end}
    \State \textbf{(Build $\mathcal L_2$)} Construct $\mathcal L_2(\mathcal P_k)$ by fixing the demand variables to $\mathbf{\tilde q^k}$
    \label{a5:l2-begin}
    \State \textbf{($\mathcal L_2$ step)} Solve $\mathcal L_2(\mathcal P_k)$ within $\mathbf{T_{max}}$ to obtain candidate $\Tilde{\mathcal X}^{(2)}_k$ (if available); let $t_k^{L_2}$ be the solve time.
    \label{a5:l2-solve-end}
    \If{$\Tilde{\mathcal X}^{(2)}_k$ is obtained}
    \label{a5:l2-recovery-begin}
        \State \textbf{(Recovery from $\mathcal L_2$)} Apply Algorithm~\ref{alg:restricted-neighborhood} with $(\mathcal N_1,\mathcal L_2(\mathcal P_k),\Tilde{\mathcal X}^{(2)}_k, \mathcal I_s,\mathbf{T_{max}})$ ;
               record the recovery time as $t_k^{rec2}$
               
        \If{Algorithm \ref{alg:restricted-neighborhood} returns a feasible solution $\mathcal S_2$ with objective $c(\mathcal S_2)>c^\star$}
            \State $\mathcal S^\star \gets \mathcal S_2$;\quad $c^\star \gets c(\mathcal S_2)$; 
            \State $t_k^{LP} \gets t_k^{L_1} + t_k^{L_2}$; $t_k^{rec} \gets t_k^{rec2}$
        \EndIf
        \EndIf
        \label{a5:l2-end}
        \If{$|c^\star - c_{\mathcal L_1(\mathcal P_k)}^\star| \le \varepsilon_{\text{opt}}$}
            \State $\mathcal S^\star$ is \text{certified optimal for} $\mathcal N_1$; 
        \EndIf
\label{a5:end-iteration}
    \State $k \gets k+1$
\EndWhile
\State \Return $\mathcal S^\star$
\end{algorithmic}
\end{algorithm}

At each refinement level $k$ of Algorithm \ref{alg:partition-refinement-tiebreak}, we first solve $\mathcal L_1(\mathcal P_k)$ to obtain an integer candidate $\Tilde{\mathcal X}^{(1)}_k$ and the associated water demand solution $\mathbf{\tilde q^k} =  \{\mathbf{q^*_{i,t}} : i \in \mathcal D,\ t \in \mathcal T\}$ (Steps \ref{a5:step-l1} - \ref{a5:step-l1-time}). 
Recovery is initiated immediately from $\Tilde{\mathcal X}^{(1)}_k$ using Algorithm \ref{alg:restricted-neighborhood} (Steps \ref{a5:l1-rec-begin} -- \ref{a5:l1-rec-end}). 
Concurrently, $\mathcal L_2(\mathcal P_k)$ is constructed by fixing its demand variables to $\mathbf{\tilde q^k}$, as described in Section \ref{subsec:milp-relaxation},
and solved within its time limit to obtain a second candidate $\Tilde{\mathcal X}^{(2)}_k$ (if available) (Steps \ref{a5:l2-begin} -- \ref{a5:l2-solve-end}). 
If such a candidate is found, a separate recovery is launched from $\Tilde{\mathcal X}^{(2)}_k$ (Step \ref{a5:l2-recovery-begin} -- \ref{a5:l2-end}). 
The incumbent solution is then updated with the best feasible solution identified at refinement level $k$ , and the process continues with the refined partition $\mathcal P_{k+1}$ until termination by optimality, time limit, or maximum refinement iterations. 


Algorithm \ref{alg:partition-refinement-tiebreak} combines partition refinement, neighborhood search, and tie-breaking to provide a robust recovery framework that systematically addresses infeasibility, computational tractability, solution quality, and robustness to solution degeneracy. 
Table~\ref{tab:alg-summary} summarizes the systematic progression from the baseline recovery method (Algorithm~1) to the complete refinement and tie-breaking framework (Algorithm~5).

\begin{table*}
\centering
\renewcommand{\arraystretch}{1.15}
\begin{tabular}{p{4.2cm} p{10.8cm}}
\toprule
\textbf{Algorithm} & \textbf{Description / Extension} \\
\midrule
\text{Alg.~1 -- Baseline Recovery} & 
Solves $\mathcal{L}_1(P)$ to obtain an integer candidate $\tilde{\mathcal{X}}$; fixes integer variables in $\mathcal{N}_1$ to $\tilde{\mathcal{X}}$ and solves the resulting NLP to recover a feasible solution. \\[4pt]

\text{Alg.~2 -- Neighborhood Search} & 
Extends Alg.~\ref{alg:baseline} by solving $\mathcal{N}_2(\tilde{X}, h)$ with increasing neighborhood size $h$ until a feasible solution is found. \\[4pt]

\text{Alg.~3 -- Restricted Neighborhood} & 
Variant of neighborhood search that restricts the Hamming-distance constraint to a subset $\mathcal{I}_s \subset \mathcal{I}$, improving computational tractability. \\[4pt]

\text{Alg.~4 -- Iterative Partition Refinement} & 
Embeds Alg.~3 within an iterative refinement framework that progressively refines partitions $\mathcal{P}_k$, recovers feasible solutions from candidates across successive partitions, and updates the incumbent solution until convergence or termination by time or iteration limit. \\[4pt]

\text{Alg.~5 -- Refinement + Tie-Breaker} & 
Extends Alg.~4 by solving $\mathcal{L}_2(\mathcal{P}_k)$ at each refinement level to resolve degeneracy and generate an additional candidate for recovery, improving robustness and consistency. \\ 
\bottomrule
\end{tabular}
\caption{Summary of the progression from Algorithm~\ref{alg:baseline} to Algorithm~\ref{alg:partition-refinement-tiebreak}. 
Each algorithm builds upon the previous one by addressing its specific limitations.}
\label{tab:alg-summary}
\end{table*}

\section{Results}
\label{sec:results}
In this section, we evaluate the effectiveness of Algorithm~\ref{alg:partition-refinement-tiebreak} in recovering feasible solutions to the demand maximization problem, formulated as $\mathcal N_1$, on five benchmark test instances from the literature. 
Table~\ref{tab:instance_summary} summarizes the instances, their sources, and key structural characteristics, including the number of junctions, pipes, pumps, tanks, and time points.

\begin{table}
\begin{tabular}{|c|c|c|cccc|c|}
\hline
\multirow{2}{*}{\#} & \multirow{2}{*}{Instance} & \multirow{2}{*}{Source}                    & \multicolumn{4}{c|}{Network Size}                                                                            & \multirow{2}{*}{$|\mathcal{T}|$} \\ \cline{4-7}
                          &                          &                                               & \multicolumn{1}{c|}{$|\mathcal{J}|$} & \multicolumn{1}{c|}{$|\mathcal{A}_p|$} & \multicolumn{1}{c|}{$|\mathcal{A}_{pu}|$} & \multicolumn{1}{c|}{$|\mathcal{TK}|$} &                                 \\ \hline
1                         & Poormond                 & \cite{tasseff2024polyhedral}
 & \multicolumn{1}{c|}{42}           & \multicolumn{1}{c|}{44}       & \multicolumn{1}{c|}{7}        & 5        & 48                              \\ \hline
2                         & Cohen                    & \cite{cohen2000optimal}                      & \multicolumn{1}{c|}{8}            & \multicolumn{1}{c|}{7}        & \multicolumn{1}{c|}{3}        & 1        & 12                              \\ \hline
3                         & VanZyl-12                   & \cite{van2004operational}                   & \multicolumn{1}{c|}{13}           & \multicolumn{1}{c|}{15}       & \multicolumn{1}{c|}{3}        & 2        & 12                               \\ \hline
4                         & VanZyl-6                   & \cite{van2004operational}                   & \multicolumn{1}{c|}{13}           & \multicolumn{1}{c|}{15}       & \multicolumn{1}{c|}{3}        & 2        & 6                              \\ \hline
5                         & AT(M)                    & \cite{tasseff2024polyhedral,rao2007use} & \multicolumn{1}{c|}{19}           & \multicolumn{1}{c|}{41}       & \multicolumn{1}{c|}{3}        & 2        & {12}                              \\ \hline
\end{tabular}
\caption{Summary of the instances used in this study}
\label{tab:instance_summary}
\end{table}

\subsection{Evaluation Metrics}
We compare Algorithm~\ref{alg:partition-refinement-tiebreak} against two benchmarks: (i) the baseline recovery algorithm (Algorithm~\ref{alg:baseline}) and (ii) spatial branch-and-bound (SBB), a global optimization method that guarantees the optimum of $\mathcal N_1$ when allowed sufficient time. 
In our experiments, Algorithm \ref{alg:baseline} failed to compute feasible solutions in most instances. Therefore, we primarily use SBB as the benchmark for performance evaluation. 
We evaluate the algorithms using two metrics: (i) the time required to compute a feasible solution and (ii) the quality of the recovered feasible solution, reported as the percentage optimality gap between the solution objective value and its upper bound. If an optimal solution is available for the instance, we use its optimal value as the upper bound; otherwise, we use the objective value of the MILP relaxation candidate (i.e., the solution of $\mathcal L_1(\cdot)$) as the upper bound.

\subsection{Implementation and Parametric Choices}
We implemented all the experiments in Julia, using SCIP~8.0.3 \cite{BestuzhevaEtal2021OO}---with the SBB method---
to solve  the MINLP formulations, $\mathcal N_1$ and $\mathcal N_2(\cdot)$, and CPLEX~22.1.0 to solve the MILP relaxations, $\mathcal L_1(\cdot)$ and $\mathcal L_2(\cdot)$. 
All computations were performed on a MacBook Pro with 16~GB RAM and a 2.8~GHz Quad-Core Intel Core i7 processor.
We imposed a time limit of $T_{\max}=3000$ seconds on each call to $\mathcal L_1(\cdot)$, $\mathcal L_2(\cdot)$, and $\mathcal N_2(\cdot)$ within the algorithm, while allowing up to 15000 seconds for solving $\mathcal N_1$ with SBB to obtain high-quality benchmarks. 



We determined $\mathcal I_s$ by trial and error and found that restricting the search to the pump activation variables and using the  explicit-direction-based formulation of $\mathcal N_1$, provided the best balance between computational efficiency and feasibility recovery in Algorithm~\ref{alg:restricted-neighborhood}. 
We therefore adopt this choice for all reported experiments. 
In this formulation, the pump activation variables represent the true controllable decisions, whereas flow-direction variables are modeling artifacts. 
Nevertheless, explicitly modeling and fixing flow-direction variables improves information transfer from the MILP relaxations to the MINLP and reduces the computational effort required to solve $\mathcal N_2(\cdot)$, thereby enhancing recovery efficiency.

\subsection{Comparison of the Best Feasible Solutions}
With the aforementioned experimental setup and parametric choices, we compute feasible solutions to $\mathcal N_1$ for each instance, using both Algorithm \ref{alg:partition-refinement-tiebreak} and the SBB method.
Then, we compare the best feasible solutions yielded by these methods by summarizing their associated metrics in Table \ref{tab:results_summary}.
For each instance, listed as a row, the table reports the problem size (number of integer variables and nonconvex constraints), the quality of the best solution (percentage gap with respect to the upper bound), the time required to compute it, and the \emph{Computation Speedup Factor}, defined as the ratio of the computation time for SBB to that for Algorithm~\ref{alg:partition-refinement-tiebreak}.

\begin{table*}
\centering
\begin{tabular}{|c|cc|cc|cc|c|}
\hline
\multirow{2}{*}{\begin{tabular}[c]{@{}c@{}}Instance\end{tabular}} & \multicolumn{2}{c|}{Size of $\mathcal N_1$} & \multicolumn{2}{c|}{\begin{tabular}[c]{@{}c@{}}Optimality Gap of the Best Solution\\ (\%gap w.r.t. Best U.B)\end{tabular}} & \multicolumn{2}{c|}{\begin{tabular}[c]{@{}c@{}}Time (seconds) to compute\\ the best solution\end{tabular}} & \multirow{2}{*}{\begin{tabular}[c]{@{}c@{}}Computation \\Speedup Factor
 \end{tabular}} \\ \cline{2-7}
 & \multicolumn{1}{c|}{$\mathbf{n_I}$} & \begin{tabular}[c]{@{}c@{}}\# {Nonconvex}\\ constraints\end{tabular} & \multicolumn{1}{c|}{SBB Method} &Algorithm \ref{alg:partition-refinement-tiebreak} & \multicolumn{1}{c|}{SBB Method} & Algorithm \ref{alg:partition-refinement-tiebreak} &  \\ \hline
Poormond & \multicolumn{1}{c|}{2448} & 2448 & \multicolumn{1}{c|}{0} & 0 & \multicolumn{1}{c|}{28.55} & 1.70 & 16.79 \\ \hline
Cohen & \multicolumn{1}{c|}{120} & 120 & \multicolumn{1}{c|}{0} & 0 & \multicolumn{1}{c|}{13.30} & 0.85 & 15.65 \\ \hline
VanZyl-12 & \multicolumn{1}{c|}{216} & 216 & \multicolumn{1}{c|}{NA} & 3.03\% & \multicolumn{1}{c|}{\textgreater{}15000} & 527.29 & \textgreater{}28.45 \\ \hline
VanZyl-6 & \multicolumn{1}{c|}{108} & 108 & \multicolumn{1}{c|}{NA} & 4.57\% & \multicolumn{1}{c|}{\textgreater{}15000} & 256.95 & \textgreater{}58.38 \\ \hline
ATM & \multicolumn{1}{c|}{528} &  528 & \multicolumn{1}{c|}{NA} & 97.17\% & \multicolumn{1}{c|}{\textgreater{}15000} & 1376.09 & \textgreater{}10.9 \\ \hline
\end{tabular}
\caption{Comparison of the best feasible solutions to $\mathcal N_1$ computed by Algorithm \ref{alg:partition-refinement-tiebreak} against the SBB method}
\label{tab:results_summary}
\end{table*}

Table~\ref{tab:results_summary} shows that both Algorithm~\ref{alg:partition-refinement-tiebreak} and SBB computed optimal solutions for the Poormond and Cohen instances. 
However, Algorithm~\ref{alg:partition-refinement-tiebreak} was faster by at least an order of magnitude, achieving speedup factors of 16.79 and 15.65, respectively (see the last column of Table~\ref{tab:results_summary}). 
For the remaining instances, SBB failed to find a single feasible solution within 15000 seconds, whereas Algorithm~\ref{alg:partition-refinement-tiebreak} determined solutions within a few hundred seconds. 
Among these, the solutions for VanZyl-6 and VanZyl-12 were near-optimal, with objective values within 5\% of the optimum. 
The solution for the ATM instance, although feasible, exhibited a large gap relative to its upper bound, indicating potential for improvement in either the solution or the bound provided by $\mathcal L_1(\cdot)$. 
Overall, Algorithm~\ref{alg:partition-refinement-tiebreak} consistently outperformed SBB by computing feasible solutions (and higher-quality solutions when applicable) at least an order of magnitude faster.


Next, we analyze the algorithm’s performance in detail by examining the iteration-wise evolution (with partition refinement) of its solutions for each instance.




\subsection{Solution Trends of Algorithm~\ref{alg:partition-refinement-tiebreak} with Partition Refinement} 
\label{subsec:trends}
Each iteration $k$ of Algorithm~\ref{alg:partition-refinement-tiebreak} corresponds to a partition $\mathcal P_k$, obtained by iteratively bisecting the initial domain $\Omega$. 
At iteration $k$, $\Omega$ is thus divided into $n_{\mathcal P_k}:= 2^k$ partitions.
Earlier, we reported the best solutions computed by the algorithm across all iterations; here, we analyze how these solutions evolve with successive refinement levels. 

To this end, for each $k$, we create the corresponding partition $\mathcal P_k$ and apply Steps~\ref{a5:step-l1}–\ref{a5:end-iteration} of Algorithm~\ref{alg:partition-refinement-tiebreak}. We summarize these results for each instance in Tables~\ref{tab:poormond}–\ref{tab:atm}. 
For each partition $\mathcal P_k$ (or iteration $k$), 
Tables~\ref{tab:poormond}–\ref{tab:atm} report the following information: 
Column~3 lists the Hamming distance ($\bar h$) between the recovered solution and the candidate. 
Columns~4–6 present the objective values of the best solutions to $\mathcal L_1(\mathcal P_k)$, Algorithm~\ref{alg:partition-refinement-tiebreak} ($c^*_{Alg. \ref{alg:partition-refinement-tiebreak}}$), and the SBB benchmark ($c^*_{SBB}$), respectively. 
Column~7 shows the optimality gap, defined as the percentage difference between the objective values of Algorithm~\ref{alg:partition-refinement-tiebreak} and $\mathcal L_1(\cdot)$. 
Columns~8–11 provide the computation times for (i) candidate generation, $t_k^{LP}$
(ii) solution recovery $t_k^{rec}$, 
(iii) their total time corresponding to Steps~\ref{a5:step-l1}–\ref{a5:end-iteration} of Algorithm~\ref{alg:partition-refinement-tiebreak}, i.e., $t_k^{LP} + t_k^{rec}$, 
and (iv) the SBB benchmark, $t_{SBB}$. 

\begin{figure}
    \centering
\includegraphics[width = 0.5\textwidth]{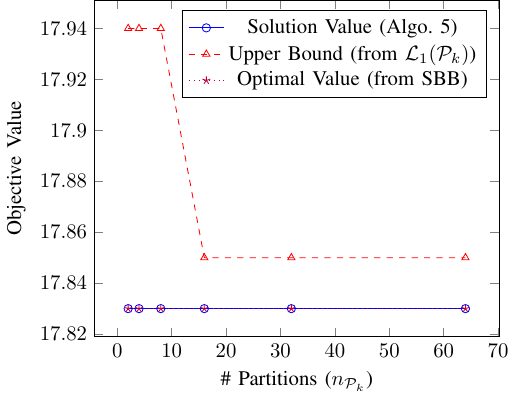}
\caption{Comparison of the solution objective values obtained by Algorithm~\ref{alg:partition-refinement-tiebreak} with the optimal value and the relaxation upper bounds as a function of the number of partitions for the Poormond instance. The algorithm achieves optimality from the first iteration and maintains negligible gaps with the upper bounds across subsequent refinements.}
    \label{fig:Poormond}
\end{figure}

Figures~\ref{fig:Poormond}, \ref{fig:cohen-12}, \ref{fig:VanZyl 6}, \ref{fig:VanZyl 12}, and \ref{fig:atm} illustrate the aforementioned results by plotting the recovered objective values together with their corresponding upper bounds as functions of the number of partitions. 
In the following subsections, we discuss the iteration-wise trends for each instance in detail.

\subsubsection{Poormond}
For the Poormond instance, Algorithm~\ref{alg:partition-refinement-tiebreak} recovered \textit{optimal solutions} from the very first iteration (with $n_{\mathcal P_k} = 2$), as indicated by the matching objective value with the SBB solution in Table~\ref{tab:poormond}. 
At each iteration, we have $\bar h = 0$, indicating that the candidate integer solution was already feasible for $\mathcal N_1$; in other words, even Algorithm~\ref{alg:baseline} was sufficient for recovery. 
The algorithm did not terminate immediately, however, because of a non-zero gap between the objective values of the feasible $\mathcal N_1$ solution and the relaxation candidate, a discrepancy attributable to the weakness of the relaxation $\mathcal L_1(\cdot)$. 
This gap was modest—within $1\%$ in the first iteration—and subsequent refinements tightened the relaxation further, reducing the gap to below $0.1\%$. 

In terms of runtime, SBB required 28.55 seconds to compute an optimal solution, whereas Algorithm~\ref{alg:partition-refinement-tiebreak} obtained the same solution in only 1.7 seconds, more than an order of magnitude faster. 
Even when initialized with finer partitions or iterated across all refinements, the total runtime of our algorithm remained lower than that of SBB despite the additional cost of candidate generation for finer partitions. 
Overall, these results show that Algorithm~\ref{alg:partition-refinement-tiebreak} matched SBB in solution quality while substantially outperforming it in speed.

\begin{table*}
\centering
\begin{tabular}{|c|c|c|ccc|c|cccc|}
\hline
\multirow{4}{*}{$k$} & \multirow{4}{*}{\begin{tabular}[c]{@{}c@{}} \# Partitions \\($n_{\mathcal P_k}$)\end{tabular}} & \multirow{4}{*}{$\bar h$}& \multicolumn{3}{c|}{Solution Objective Value} & \multirow{4}{*}{\begin{tabular}[c]{@{}c@{}}  Optimality Gap\\($c^*_{Alg. \ref{alg:partition-refinement-tiebreak}}$ vs  $c_{\mathcal L_1(\mathcal P_k)}^\star$)
\end{tabular}} & \multicolumn{4}{c|}{Computation Time (seconds)}  \\ \cline{4-6}\cline{8-11}
&  & & $c_{\mathcal L_1(\mathcal P_k)}^\star$ & $c^*_{Alg.\ref{alg:partition-refinement-tiebreak}}$ & $c^*_{SBB}$ &  & \begin{tabular}[c]{@{}c@{}} Candidate \\Generation \\ ($t_k^{LP}$)\end{tabular}& \begin{tabular}[c]{@{}c@{}} Solution \\Recovery\\($t_k^{rec}$)\end{tabular}  & \begin{tabular}[c]{@{}c@{}}Total\\($t_k^{LP} + t_k^{rec}$)\end{tabular} & \begin{tabular}[c]{@{}c@{}}$t_{SBB}$\end{tabular} \\ \hline
1 & 2 & 0  & 17.94 & 17.83 &\multirow{6}{*}{17.83}& 0.61 \% & 0.28 & 1.42 & 1.70 & \multirow{6}{*}{28.55} \\
2 &4 &0 &17.94 & 17.83   & & 0.61 \% & 0.48 & 1.52 & 2.00 &  \\
3 &8 & 0 & 17.94 & 17.83   & & 0.61 \% & 0.78 & 1.39 & 2.17 &  \\
4 &16 & 0  & 17.85 & 17.83  &  & 0.08 \% & 1.45 & 1.41 & 2.86 &  \\
5 &32 & 0  & 17.85 & 17.83  & & 0.08 \% & 3.37 & 1.36 & 4.73 &  \\
6 &64 & 0  & 17.85 & 17.83  &  & 0.08 \% & 7.90 & 1.34 & 9.24 &  \\ \hline
\end{tabular}
\caption{
Iteration-wise results for the Poormond instance with table structure as described above in Section~\ref{subsec:trends}.}
\label{tab:poormond}
\end{table*}

\subsubsection{Cohen}
\begin{figure}
\includegraphics[width=0.49\textwidth]{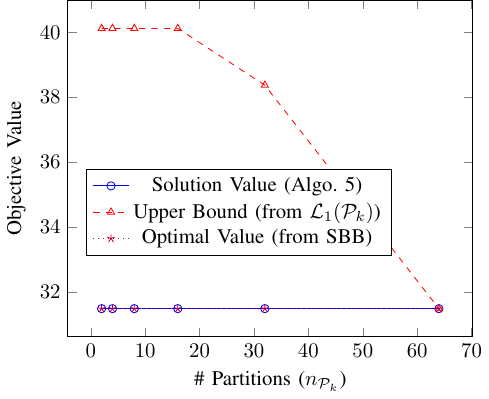}
   \caption{Comparison of the objective values obtained by Algorithm~\ref{alg:partition-refinement-tiebreak} with the optimal value and the relaxation upper bounds as a function of the number of partitions for the Cohen instance. The solutions are optimal across all refinement levels, while the upper bound converges progressively toward the optimal value with increasing partition granularity.}

    \label{fig:cohen-12}
\end{figure}

For the Cohen instance, Algorithm~\ref{alg:partition-refinement-tiebreak} recovered optimal solutions at all iterations 
, as shown in Figure~\ref{fig:cohen-12} and Table~\ref{tab:cohen}. 
Unlike the Poormond instance, the relaxation upper bound converged gradually to the recovered solution value with successive refinements. 
For coarser partitions ($k \leq 5$), $\bar h$ was non-zero, ranging from 12 to 9 (see column 3 of Table~\ref{tab:cohen}), indicating that a neighborhood search was required for recovery. 
This demonstrates that the baseline method (Algorithm~\ref{alg:baseline}) fails to provide feasible solutions, whereas Algorithm~\ref{alg:partition-refinement-tiebreak} consistently recovered feasible solutions at all refinement levels, many of which were optimal.

In terms of runtime, the fastest recovery of an optimal solution occurred at $k=3$ (eight partitions), requiring only 0.85 seconds—15.65 times faster than SBB. 
Even at $k=2,5,$ and $6$, Algorithm~\ref{alg:partition-refinement-tiebreak} was faster than SBB; the latter failed to yield feasible solutions when allowed the same computation times as the former. 
At $k=1$ and $k=4$, the runtime exceeded SBB because SCIP incurred overhead in certifying the optimality or infeasibility of $\mathcal N_2(\cdot)$ for certain values of $h$.

Overall, for the Cohen instance, Algorithm~\ref{alg:partition-refinement-tiebreak} outperformed Algorithm~\ref{alg:baseline} by consistently recovering optimal solutions and outperformed SBB in computation time at several refinement levels. 
Notably, Cohen is the only other instance for which SBB produced feasible solutions; for all subsequent instances, SBB failed to return a feasible solution within the 15,000-second time limit.

\begin{table*}
\centering
\begin{tabular}{|c|c|c|ccc|c|cccc|}
\hline
\multirow{4}{*}{$k$} & \multirow{4}{*}{\begin{tabular}[c]{@{}c@{}} \# Partitions \\($n_{\mathcal P_k}$)\end{tabular}} & \multirow{4}{*}{$\bar h$}& \multicolumn{3}{c|}{Solution Objective Value} & \multirow{4}{*}{\begin{tabular}[c]{@{}c@{}}  Optimality Gap\\($c^*_{Alg. \ref{alg:partition-refinement-tiebreak}}$ vs  $c_{\mathcal L_1(\mathcal P_k)}^\star$)
\end{tabular}} & \multicolumn{4}{c|}{Computation Time (seconds)}  \\ \cline{4-6}\cline{8-11}
&  & & $c_{\mathcal L_1(\mathcal P_k)}^\star$ & $c^*_{Alg.\ref{alg:partition-refinement-tiebreak}}$ & $c^*_{SBB}$ &  & \begin{tabular}[c]{@{}c@{}} Candidate \\Generation \\ ($t_k^{LP}$)\end{tabular}& \begin{tabular}[c]{@{}c@{}} Solution \\Recovery\\($t_k^{rec}$)\end{tabular}  & \begin{tabular}[c]{@{}c@{}}Total\\($t_k^{LP} + t_k^{rec}$)\end{tabular} & \begin{tabular}[c]{@{}c@{}}$t_{SBB}$\end{tabular} \\ \hline
1 & 2 & 12 & 40.12 & 31.48 & \multirow{6}{*}{31.48} & 21.54 & 0.02 & 18.30 & 18.32 & \multirow{6}{*}{13.30}\\
2 & 4 & 12 & 40.12 & 30.48 & & 24.46 & 0.17 & 0.56 & 0.73 &  \\
3 & 8 & 12 & 40.12 & 31.48 & & 21.54 & 0.16 & 0.69 & 0.85 &  \\
4 & 16 & 12 & 40.12 & 31.48 & & 21.54 & 0.33 & 18.76 & 19.09 &  \\
5 & 32 & 9 & 38.37 & 31.48 & & 17.95 & 0.57 & 8.92 & 9.49 &  \\
6 & 64 & 0 & 31.48 & 31.48 & & 0 & 1.85 & 0.22 & 2.07 &  \\ \hline
\end{tabular}
\caption{Iteration-wise results for the Cohen instance with table structure as described above in Section~\ref{subsec:trends}.}
\label{tab:cohen}
\end{table*}

\subsubsection{VanZyl Instances}
\begin{figure}
\includegraphics[width=0.49\textwidth]{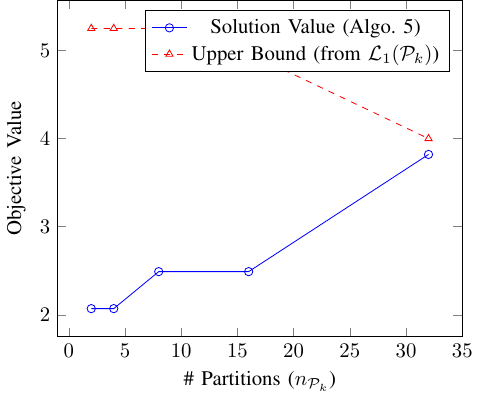}
\caption{Comparison of the objective values obtained by Algorithm~\ref{alg:partition-refinement-tiebreak} with the relaxation upper bound as a function of the number of partitions for the VanZyl-6 instance. Both the recovered solution values and the relaxation upper bounds improve with partition refinement, reducing the optimality gap to below 5\% and demonstrating monotonic convergence behavior, whereas the SBB method fails to compute a feasible solution within the allowed time.}

    \label{fig:VanZyl 6}
\end{figure}

\begin{figure}
\includegraphics[width=0.5\textwidth]{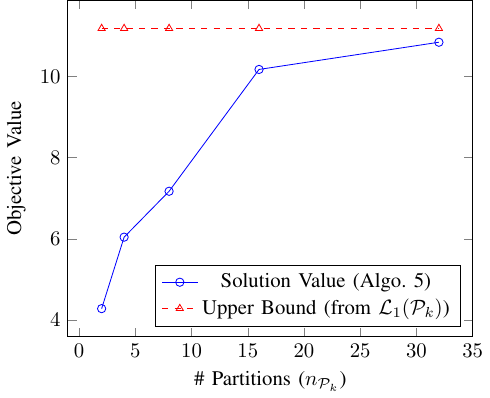}
\caption{Comparison of the objective values obtained by Algorithm~\ref{alg:partition-refinement-tiebreak} with the relaxation upper bound as a function of the number of partitions for the VanZyl-12 instance. While the relaxation upper bound remains constant, the recovered solution value improves progressively with partition refinement, closing the optimality gap to approximately~3\%, whereas the SBB method fails to compute a feasible solution within the allowed time.}
    \label{fig:VanZyl 12}
\end{figure}

The two VanZyl instances are derived from the same network but differ in planning horizon and demand profile patterns. 
For both instances, SBB failed to find feasible solutions within the 15,000-second time limit, whereas Algorithm~\ref{alg:partition-refinement-tiebreak} consistently recovered feasible solutions for all partition levels. 
The recovered solution values improved monotonically with partition refinement, progressively closing the gap with their upper bounds, as shown in Figures~\ref{fig:VanZyl 6} and~\ref{fig:VanZyl 12}. 
For VanZyl-6, the upper bounds initially remained unchanged with refinement but gradually improved, reducing the final gap to 4.57\%. 
In contrast, for VanZyl-12, the upper bounds remained constant throughout, and the gap closed to 3.03\% solely through improvements in solution quality.

As shown in Tables~\ref{tab:vanzyl6} and~\ref{tab:vanzyl12}, similar to the Cohen instance, $\bar h \ne 0$ for $k < 5$, indicating that the baseline recovery algorithm failed to yield feasible solutions with coarse partitioning, thus necessitating the neighborhood search in Algorithm~\ref{alg:partition-refinement-tiebreak}. 

The computation times of Algorithm~\ref{alg:partition-refinement-tiebreak} ranged from a few seconds to several hundred seconds, yet it consistently outperformed SBB, which was unable to find feasible solutions within the allotted time. 
These results demonstrate the algorithm’s robustness and effectiveness in rapidly recovering feasible solutions for complex networks.


\begin{table*}
\centering
\begin{tabular}{|c|c|c|ccc|c|cccc|}
\hline
\multirow{4}{*}{$k$} & \multirow{4}{*}{\begin{tabular}[c]{@{}c@{}} \# Partitions \\($n_{\mathcal P_k}$)\end{tabular}} & \multirow{4}{*}{$\bar h$}& \multicolumn{3}{c|}{Solution Objective Value} & \multirow{4}{*}{\begin{tabular}[c]{@{}c@{}}  Optimality Gap\\($c^*_{Alg. \ref{alg:partition-refinement-tiebreak}}$ vs  $c_{\mathcal L_1(\mathcal P_k)}^\star$)
\end{tabular}} & \multicolumn{4}{c|}{Computation Time (seconds)}  \\ \cline{4-6}\cline{8-11}
&  & & $c_{\mathcal L_1(\mathcal P_k)}^\star$ & $c^*_{Alg.\ref{alg:partition-refinement-tiebreak}}$ & $c^*_{SBB}$ &  & \begin{tabular}[c]{@{}c@{}} Candidate \\Generation \\ ($t_k^{LP}$)\end{tabular}& \begin{tabular}[c]{@{}c@{}} Solution \\Recovery\\($t_k^{rec}$)\end{tabular}  & \begin{tabular}[c]{@{}c@{}}Total\\($t_k^{LP} + t_k^{rec}$)\end{tabular} & \begin{tabular}[c]{@{}c@{}}$t_{SBB}$\end{tabular} \\ \hline
1 & 2 & 2 & 5.25 & 2.07 & \multirow{6}{*}{NA} & 60.55 & 0.24 & 1.18 & 1.42 & \multirow{5}{*}{$>15000$}\\
2 & 4 & 2 & 5.25 & 2.07 & & 60.5 & 0.35 & 5.02 & 5.37 & \\
3 & 8 & 3 & 5.25 & 2.49 & & 52.58 & 1.10 & 11.01 & 12.11 & \\
4 & 16 & 3 & 4.97 & 2.49 & & 49.97 & 31.27 & 8.37 & 39.64 &  \\
5 & 32 & 0 & 4.0 & 3.82 & & 4.57 & 256.80 & 0.15 & 256.95 &  \\
\hline
\end{tabular}
\caption{Iteration-wise results for the VanZyl-6 instance with table structure as described above in Section~\ref{subsec:trends}.}
\label{tab:vanzyl6}
\end{table*}

\begin{table*}
\centering
\begin{tabular}{|c|c|c|ccc|c|cccc|}
\hline
\multirow{4}{*}{$k$} & \multirow{4}{*}{\begin{tabular}[c]{@{}c@{}} \# Partitions \\($n_{\mathcal P_k}$)\end{tabular}} & \multirow{4}{*}{$\bar h$}& \multicolumn{3}{c|}{Solution Objective Value} & \multirow{4}{*}{\begin{tabular}[c]{@{}c@{}}  Optimality Gap\\($c^*_{Alg. \ref{alg:partition-refinement-tiebreak}}$ vs  $c_{\mathcal L_1(\mathcal P_k)}^\star$)
\end{tabular}} & \multicolumn{4}{c|}{Computation Time (seconds)}  \\ \cline{4-6}\cline{8-11}
&  & & $c_{\mathcal L_1(\mathcal P_k)}^\star$ & $c^*_{Alg.\ref{alg:partition-refinement-tiebreak}}$ & $c^*_{SBB}$ &  & \begin{tabular}[c]{@{}c@{}} Candidate \\Generation \\ ($t_k^{LP}$)\end{tabular}& \begin{tabular}[c]{@{}c@{}} Solution \\Recovery\\($t_k^{rec}$)\end{tabular}  & \begin{tabular}[c]{@{}c@{}}Total\\($t_k^{LP} + t_k^{rec}$)\end{tabular} & \begin{tabular}[c]{@{}c@{}}$t_{SBB}$\end{tabular} \\ \hline
1 & 2 & 4 &  11.18 & 4.28 & \multirow{6}{*}{NA} & 61.69 & 56.75 &  5.11& 61.86 & \multirow{5}{*}{$> 15000$}\\
2 & 4 & 3 & 11.18 & 6.04 & & 45.99 & 642.09 & 39.79& 681.88 &  \\
3 & 8 & 8 & 11.18 & 7.17 & & 36.49 & 93.17 & 28.31 & 121.48 &  \\
4 & 16 & 1 & 11.18 & 10.18 & & 8.92 & 27.74 & 18.64 & 46.38 &  \\
5 & 32 & 0 & 11.18 & 10.84 & & 3.03 & 414.11 & 13.18 & 527.29 &  \\
\hline
\end{tabular}
\caption{Iteration-wise results for the VanZyl-12 instance with table structure as described above in Section~\ref{subsec:trends}.}
\label{tab:vanzyl12}
\end{table*}


\subsubsection{ATM}
\begin{figure}
\includegraphics[width = 0.49\textwidth]{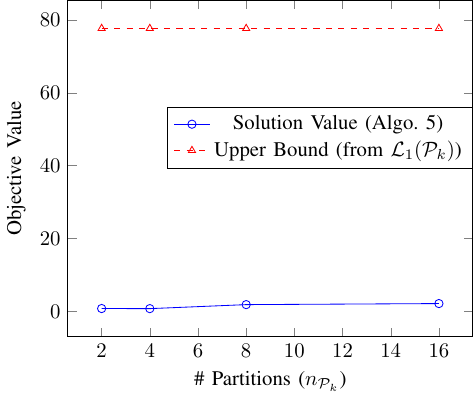}
\caption{Comparison of the objective values obtained by Algorithm~\ref{alg:partition-refinement-tiebreak} with the relaxation upper bounds as a function of the number of partitions for the ATM instance. Although the gaps remain large across refinements, Algorithm~\ref{alg:partition-refinement-tiebreak} successfully recovers feasible solutions at all partition levels, whereas the SBB method fails to obtain an optimal solution within the allowed time.}
    \label{fig:atm}
\end{figure}

Similar to the VanZyl instances, SBB failed to compute feasible solutions for the ATM instance, whereas Algorithm~\ref{alg:partition-refinement-tiebreak} successfully recovered feasible solutions for all partition levels. 
The recovered solution values and their gaps with respect to the upper bounds are presented in Figure~\ref{fig:atm} and Table~\ref{tab:atm}. 
As shown in the table, $\bar h \ne 0$, consistent with the previous instances, indicating the failure of the baseline recovery algorithm and the necessity of the neighborhood search in Algorithm~\ref{alg:partition-refinement-tiebreak}.

Although all recovered solutions were feasible, the gaps with respect to the upper bounds remained large compared to the other instances (Figure~\ref{fig:atm}). 
Even at $k = 4$, the gap reached only 97.17\%. 
This large gap may be attributed to weak upper bounds in the MILP relaxations or to lower-quality feasible solutions in the nonlinear subproblems.

\begin{table*}
\centering
\begin{tabular}{|c|c|c|ccc|c|cccc|}
\hline
\multirow{4}{*}{$k$} & \multirow{4}{*}{\begin{tabular}[c]{@{}c@{}} \# Partitions \\($n_{\mathcal P_k}$)\end{tabular}} & \multirow{4}{*}{$\bar h$}& \multicolumn{3}{c|}{Solution Objective Value} & \multirow{4}{*}{\begin{tabular}[c]{@{}c@{}}  Optimality Gap\\($c^*_{Alg. \ref{alg:partition-refinement-tiebreak}}$ vs  $c_{\mathcal L_1(\mathcal P_k)}^\star$)
\end{tabular}} & \multicolumn{4}{c|}{Computation Time (seconds)}  \\ \cline{4-6}\cline{8-11}
&  & & $c_{\mathcal L_1(\mathcal P_k)}^\star$ & $c^*_{Alg.\ref{alg:partition-refinement-tiebreak}}$ & $c^*_{SBB}$ &  & \begin{tabular}[c]{@{}c@{}} Candidate \\Generation \\ ($t_k^{LP}$)\end{tabular}& \begin{tabular}[c]{@{}c@{}} Solution \\Recovery\\($t_k^{rec}$)\end{tabular}  & \begin{tabular}[c]{@{}c@{}}Total\\($t_k^{LP} + t_k^{rec}$)\end{tabular} & \begin{tabular}[c]{@{}c@{}}$t_{SBB}$\end{tabular} \\ \hline
1 & 2 & 36& 77.63 & 0.84 & \multirow{4}{*}{NA} & 98.92 & 0.08 & 5.71 & 5.79 & \multirow{4}{*}{$> 15000$} \\
2 &4 & 9 & 77.63 & 0.81 & & 98.95 & 2.30 & 12.10 & 14.40 &  \\
3 &8 & 10 & 77.63 & 1.91 & & 97.45 & 152.68 & 3000 ( = $\mathbf{T_{max}}$) & {152.68 + $\mathbf{T_{max}}$} &  \\
4 &16 & 8 & 77.63 & 2.20 & & 97.17 & 1366.75 & 9.34 & 1376.09 &   \\
\hline
\end{tabular}
\caption{Iteration-wise results for the ATM instance with table structure as described above in Section~\ref{subsec:trends}.}
\label{tab:atm}
\end{table*}


Across all instances, Algorithm~\ref{alg:partition-refinement-tiebreak} consistently recovered feasible solutions for every partition level, even in cases where both the baseline recovery method and SBB failed. 
For smaller networks (Poormond and Cohen), the algorithm achieved optimal solutions—often at the coarsest partitions—and was more than an order of magnitude faster than SBB. 
For the larger and more complex VanZyl and ATM instances, SBB failed to produce feasible solutions within the time limit, whereas Algorithm~\ref{alg:partition-refinement-tiebreak} successfully computed feasible solutions at all refinement levels. 
In these cases, solution quality improved monotonically with partition refinement, progressively tightening the gaps with respect to the relaxation bounds. 
Overall, the results highlight the robustness and computational efficiency of the proposed algorithm across problem scales. 


\section{Conclusions}
\label{sec:summary}
This study presented a robust algorithm for recovering feasible solutions to nonconvex MINLPs from their MILP relaxations. The algorithm was motivated by the demand maximization problem in water distribution networks (WDNs), where nonlinear hydraulic physics and discrete operational decisions result in large-scale nonconvex MINLPs that are challenging to solve using existing global optimization techniques.

The proposed algorithm integrates three complementary components---partition refinement, neighborhood search, and a tie-breaking strategy---to systematically address infeasibility, relaxation weakness, tractability, and solution degeneracy. Partition refinement progressively strengthens MILP relaxations while maintaining tractability; neighborhood search repairs infeasible candidates from coarse relaxations by exploiting integer-variable proximity; and the embedded tie-breaking MILP ensures consistent candidate selection and enables additional recoveries. Together, these elements form a unified framework that efficiently converts MILP relaxations of any granularity into high-quality feasible solutions.

Extensive numerical experiments on benchmark WDNs demonstrate that the algorithm reliably recovers feasible and often near-optimal solutions across diverse network sizes and topologies, even when baseline methods fail. When global methods such as spatial branch-and-bound (SBB) are able to compute optimal solutions, the proposed method computes them more than an order of magnitude faster. For larger and more complex instances, where SBB fails to obtain a feasible solution within the allowed time, the algorithm consistently identifies feasible solutions with small optimality gaps. The results also corroborate the monotonic improvement of solution quality with progressive partition refinement.

The framework is modular and broadly applicable: it can operate as a standalone recovery procedure or as a subroutine in decomposition and global optimization schemes for large-scale MINLPs.
Beyond water systems, the framework extends naturally to other infrastructure optimization problems governed by nonlinear physics and discrete decisions, including natural gas, power, and district heating networks.
Future work will focus on improving scalability through decomposition and parallelization, strengthening MILP relaxations via adaptive bound tightening, and accelerating the integer neighborhood search.

Overall, the proposed algorithm offers a practical and computationally efficient pathway for obtaining feasible, high-quality solutions to large-scale nonconvex MINLPs where exact global optimization remains computationally prohibitive.

\section*{Acknowledgment}
The authors acknowledge the support provided by the following two projects to accomplish this work.
1) LANL's Directed Research \& Development funded project, ``20240564ECR: Scalable Algorithms for Solving Mixed Integer Nonlinear Programs With a Separable Structure", and 
2) Department of Energy's Advanced Grid Modeling project, ``Coordinated planning and operation of water and power infrastructures for increased resilience and reliability".
The research conducted at Los Alamos National Laboratory is done under the auspices of the National Nuclear Security Administration of the U.S. Department of Energy under Contract No. 89233218CNA000001.
This article has been approved for unlimited public release with the report number LA-UR-25-20735.
\bibliographystyle{IEEEtran}
\bibliography{references.bib}

\appendix

\input{notation}

\end{document}

%% file: notation.tex
\subsection{Lists of Definitions}
For the convenience of readers, here we summarize the list of all decision variables, parameters, sets, and functions used in Section \ref{subsec:water-flow-models} to describe the mathematical models for water flow through distribution networks 

\subsubsection{Sets \& Functions}
\noindent\\
$\mathcal{J}$ - set of all junctions in the network\\
$\mathcal{A}_p$ - set of pipes in the network\\
$\mathcal{A}_{pu}$ - set of pumps in the network\\
$\mathcal{A}$ - set of all arcs in the networks; $\mathcal{A}$ = $\mathcal{A}_p \cup \mathcal{A}_{pu}$\\
$\mathcal{D}$ - set of all demand points in the network\\
$\mathcal{TK}$ - set of all tanks in the network\\
$\mathcal{R}$ - set of all reservoirs in the network\\
$\mathcal{T}$ - set of all time points in the planning horizon\\
$to(a)$ - \textit{`to'} junction of arc $a \in \mathcal{A}$\\
$fr(a)$ - \textit{`from'} junction of arc $a \in \mathcal{A}$\\
$inarcs(j)$ - set of arcs for which $j \in \mathcal{J}$ is the `to' junction\\
$outarcs(j)$ - set of arcs for which $j \in \mathcal{J}$ is the `from' junction\\
$tk(i)$ - set of tanks connected to junction $j \in \mathcal{J}$\\
$r(i)$ - set of reservoirs connected to junction $j \in \mathcal{J}$\\
$d(i)$ - set of demand points connected to junction $j \in \mathcal{J}$

\subsubsection{Parameters}
\noindent\\
$\mathbf{\overline{h}_j}$ - upper bound on the head at junction $j \in \mathcal{J}$ \\
$\mathbf{\underline{h}_j}$ - upper bound on the head at junction $j \in \mathcal{J}$ \\
$\mathbf{\overline{q}_{a}^+}$ - upper bound on the flow rate from $fr(a)$ to $to(a)$ for pipe $a \in \mathcal{A}_p$\\
$\mathbf{\overline{q}_{a}^-}$ - upper bound on the flow rate from $to(a)$ to $fr(a)$ for pipe $a \in \mathcal{A}_p$\\
$ \mathbf{\underline{q}_{a}^+}$ - lower bound on the flow rate from $fr(a)$ to $to(a)$ for pipe $a \in \mathcal{A}_p$\\
$\mathbf{\underline{q}_{a}^-}$ - lower bound on the flow rate from $to(a)$ to $fr(a)$ for pipe $a \in \mathcal{A}_p$\\
$\mathbf{\overline{\Delta h}_{a}^+}$ - upper bound on the difference of head at $to(a)$ from that at $fr(a)$ of arc a $a \in \mathcal{A}$ \\
$\mathbf{\overline{\Delta h}_{a}^-}$ - upper bound on the difference of head at $fr(a)$ from that at $to(a)$  of arc a $a \in \mathcal{A}$ \\
$\mathbf{\underline{\Delta h}_{a}^-}$ - lower bound on the difference of head at $fr(a)$ from that at $to(a)$ of arc a $a \in \mathcal{A}$ \\
$\mathbf{L_a}$ - length of pipe $a \in \mathcal{A}_p$\\
$\mathbf{r_a}$ - resistance per unit length of pipe $a \in \mathcal{A}_p$\\
$\mathbf{\overline{q}_{a}}$ - upper bound on the flow rate through pump $a \in \mathcal{A}_{pu}$ \\
$\mathbf{\underline{q}_{a}}$ - lower bound on the flow rate through pump $a \in \mathcal{A}_{pu}$\\
$\boldsymbol{\alpha_a}$, $\boldsymbol{\beta_a}$, $\boldsymbol{\gamma_a}$ - coefficients determined from the head curve of pump $a \in \mathcal{A}_{pu}$\\
$\boldsymbol{\omega_a}$, $\boldsymbol{\mu_a}$- coefficients determined from the power curve of pump $a \in \mathcal{A}_{pu}$\\
$\boldsymbol{\$_{a,t}}$ - energy price for pump $a$ at time point $t$\\
$\mathbf{A_i}$ - cross-sectional area of tank $i \in \mathcal{TK}$\\
$\mathbf{b_i}$ - elevation of the bottom of tank $i \in \mathcal{TK}$\\
$\mathbf{\overline{V}_{i,0}}$ - initial volume of water in tank $i \in \mathcal{TK}$\\
$\mathbf{\overline{V}_{i}}$ - upper bound on the volume of water in tank $i \in \mathcal{TK}$\\
$\mathbf{\underline{V}_{i}}$ - lower bound on the volume of water in tank $i \in \mathcal{TK}$\\
$\mathbf{T_f}$ - final time point in the set $\mathcal{T}$\\
$\mathbf{{d}_{i,t}}$ - maximum required withdrawal rate at demand point $i \in \mathcal{D}$ at time point $t \in \mathcal{T}$

\subsubsection{Variables}
\noindent\\
$h_{j,t}$ - the total hydraulic heat at junction $j \in \mathcal{J}$ at time point $t \in \mathcal{T}$\\
$y_{a,t}$ - binary decision variable indicating the water flow direction along pipe $a \in \mathcal{A}_p$ at time point $t \in \mathcal{T}$; $y_{a,t}$ = 1 if water flows from $fr(a)$ to $to(a)$ at time point $t$ and $y_{a,t}$ = 0 otherwise\\
$q_{a,t}$ - volumetric water flow rate through arc $a \in \mathcal{A}$ from $fr(a)$ to $to(a)$ at time point $t \in \mathcal{T}$\\
$q_{a,t}^+$ - non-negative water flow rate through pipe $a \in \mathcal{A}_p$ from $fr(a)$ to $to(a)$ at time point $t \in \mathcal{T}$\\
$q_{a,t}^-$ - non-negative water flow rate through pipe $a \in \mathcal{A}_p$ from $to(a)$ to $fr(a)$ at time point $t \in \mathcal{T}$\\
$\Delta h_{a,t}^+ $ - non-negative head difference of junction $to(a)$ from junction $fr(a)$ at time point $t \in \mathcal{T}$\\
$\Delta h_{a,t}^- $ - non-negative head difference of junction $fr(a)$ from junction $to(a)$ at time point $t \in \mathcal{T}$\\
$f_{d,t}$ - water withdrawal rate at demand point $d \in \mathcal{D}$ at time point $t \in \mathcal{T}$\\
$z_{a,t}$ - binary variable indicating the activation status of pump $a \in \mathcal{A}_{pu}$ at time point $t \in \mathcal{T}$; it takes the value 1 if the pump is active at time point $t$ and 0 otherwise\\
$g_{a,t}$ - head gain offered by pump $a \in \mathcal{A}_{pu}$ at time point $t \in \mathcal{T}$\\
$P_{a,t}$ - power consumed by pump $a \in \mathcal{A}_{pu}$ at time point $t \in \mathcal{T}$\\
$V_{i,t}$ - volume of water in tank  $i \in \mathcal{TK}$ at time point $t \in \mathcal{T}$\\
$q_{i,t}$ - water flow rate from tank $i \in \mathcal{TK}$ to the network at time point $t \in \mathcal{T}$\\
$q_{i,t}$ - water flow rate from reservoir $i \in \mathcal{R}$ to the network at time point $t \in \mathcal{T}$\\
$q_{i,t}$ - water flow rate from the network to the demand point $i \in \mathcal{D}$ to the network at time point $t \in \mathcal{T}$